\documentclass[11pt,a4paper]{article}
\usepackage{fullpage}

\usepackage{microtype}
\usepackage{graphicx}

\usepackage[colorlinks]{hyperref}
\hypersetup{
    linkcolor=red,
    citecolor=teal,
    filecolor=black,
    urlcolor=black,
}
\newcommand{\Soft}{\text{Soft}}
\newcommand{\sign}{\text{sign}}
\renewcommand{\Re}{\operatorname{Re}}
\renewcommand{\Im}{\operatorname{Im}}
\usepackage{verbatim}
\usepackage{amsthm,amsmath}
\usepackage{amssymb}
\usepackage{bm}
\usepackage{bbm}
\usepackage{dsfont}
\usepackage{url}
\usepackage{comment}
\usepackage{multirow}
\usepackage{color}
\usepackage{tikz}
\usepackage{float}
\usepackage{mathtools}
\usepackage{pifont}
\usepackage{amsmath}
\usepackage{enumerate, algorithm,algorithmic,caption,graphicx,color}
\usepackage{chngcntr}

\newcommand{\half}{ \mbox{\small$\frac{1}{2}$}}

\def\Res{{\hbox{\textup{Res}}}}

\newcommand{\be}{\begin{eqnarray}}
\newcommand{\ee}[1]{\label{#1}\end{eqnarray}}
\newcommand{\nn}{\nonumber \\}
\newcommand{\ese}{\end{eqnarray*}}
\newcommand{\bse}{\begin{eqnarray*}}

\def\qed{\hfill$\square$}

\def\cO{O}

\def\X{{\mathcal{X}}}

\def\argmin{\mathop{\hbox{\rm argmin}}}

\def\E{{\mathbf{E}}}
\def\cA{{\cal A}}
\def\N{{\mathcal{N}}}

\def\S{{\mathcal{S}}}

\def\C{\mathds{C}}
\def\R{\mathds{R}}
\def\Z{\mathds{Z}}

\newcommand{\wh}[1]{{\widehat{#1}}}

\def\Argmin{\mathop{\hbox{\rm Argmin}}}

\def\prox{\textup{Prox}}

\def\T{{\textup{T}}}
\def\H{{\textup{H}}}
\def\vphi{\varphi}

\newcommand{\cN}{\mathcal{N}}

\newcommand{\ie}{\textit{i.e.}}
\newcommand{\eg}{\textit{e.g.}}

\newcommand{\dimacorr}[2]{{\color{magenta}#2}}

\newtheorem{theorem}{Theorem}[section]

\newtheorem{lemma}{Lemma}[section]

\newtheorem{prop}{Proposition}[section]

\newcommand\DualityGap{\textup{DualityGap}}

\renewcommand{\algorithmicrequire}{\textbf{Input:}}
\renewcommand{\algorithmicensure}{\textbf{Output:}}

\newcommand{\veps}{\varepsilon}

\newcommand{\A}{\cA}

\newcommand{\fast}{\textup{fast}}
\renewcommand{\Vec}{\textup{Vec}}

\newcommand{\SNR}{\textup{SNR}}
\newcommand{\PSNR}{\textup{PSNR}}
\newcommand\RanSin{\textit{Random}}
\newcommand\CohSin{\textit{Coherent}}
\newcommand\ModSin{\textit{Modulated}}
\newcommand\fine{\textup{fine}}
\newcommand\coarse{\textup{coarse}}
\allowdisplaybreaks

\begin{document}

\title{Efficient First-Order Algorithms for Adaptive Signal Denoising}
\author{
	Dmitrii Ostrovskii \\
	INRIA -- SIERRA project-team \\
	\texttt{dmitrii.ostrovskii@inria.fr}
	\and
	Zaid Harchaoui\\
	University of Washington\\
	\texttt{zaid@uw.edu}
}
\maketitle

\begin{abstract}

We consider the problem of discrete-time signal denoising,
focusing on a specific family of non-linear convolution-type estimators.
Each such estimator is associated with a time-invariant filter which is obtained adaptively, by solving a certain convex optimization problem.
Adaptive convolution-type estimators were demonstrated to have favorable statistical properties, see~\cite{jn1-2009,jn2-2010,harchaoui2015adaptive,ostrovsky2016structure}.
Our first contribution is an efficient algorithmic implementation of these estimators via the known first-order proximal algorithms. 
Our second contribution is a computational complexity analysis of the proposed procedures, which takes into account their statistical nature and the related notion of statistical accuracy.
The proposed procedures and their analysis are illustrated on a simulated data benchmark. 

\end{abstract}
\section{Introduction} 
\label{sec:intro}

We consider the problem of discrete-time signal denoising.
The goal is to estimate a discrete-time complex signal $(x_\tau)$ observed in complex Gaussian noise of level~$\sigma$ on $[-n,n]$:
\begin{equation}
\label{eq:intro-observations}
y_\tau := x_\tau + \sigma \zeta_\tau, \quad \tau = -n, ..., n.
\end{equation}
Here, $\zeta_\tau$ are i.i.d. random variables with standard complex Gaussian distribution $\C\cN(0,1)$, that is, $\Re(\zeta_\tau)$ and $\Im(\zeta_\tau)$ are independent standard Gaussian random variables. 

Signal denoising is a classical problem in statistical estimation and signal processing; see~\cite{ik1981,nemirovski_topics,tsybakov_mono,wasserman2006all,haykin1991adaptive,kay1993fundamentals}. 
The conventional approach is to assume that $x$ comes from a known set $\X$ with a simple structure that can be exploited
to build the estimator. For example, one might consider signals belonging to linear subspaces $\S$ of signals whose spectral representation, as given by the Discrete Fourier or Discrete Wavelet transform, comes from a linearly transformed $\ell_p$-ball, see~\cite{tsybakov_mono,johnstone-book}. 
In all these cases, estimators with near-optimal statistical performance can be found in explicit form, and correspond to linear functionals of the observations $y$ -- hence the name \textit{linear estimators}.

We focus here on a family of \emph{non-linear} estimators with larger applicability and strong theoretical guarantees, in particular when the structure of the signal is unknown beforehand, as studied in~\cite{nemirovsky1991nonparametric,jn1-2009,jn2-2010,harchaoui2015adaptive,ostrovsky2016structure}.
Assuming for convenience that one must estimate~$x_t$ on~$[0, n]$ from observations~\eqref{eq:intro-observations}, these estimators can be expressed as
\begin{equation}
\label{eq:intro-convolution-est}
\wh x^\vphi_t = [\vphi * y]_t := \sum_{\tau \in \Z} \vphi_\tau y_{t-\tau} \quad 0 \le t \le n;
\end{equation}
here $\vphi$ is called a \textit{filter} and is supported on $[0,n]$ which we write as $\vphi \in \C_n(\Z)$, and $*$ is the (non-circular) discrete convolution.
For estimators in this family, the filter is then obtained as an optimal solution to some convex optimization problem.
For instance, the \textit{Penalized Least-Squares} estimator~\cite{ostrovsky2016structure} is defined by
\begin{equation}
\label{opt:l2pen-0}
\wh \vphi \in \Argmin\limits_{\vphi \in \C_n(\Z)}
\half \|F_n[y-\vphi*y]\|_{2}^2 + \lambda\|F_n[\vphi]\|_{1},
\end{equation}
where $F_n$ is the Discrete Fourier transform (DFT) on $\C^{n+1}$, and $\|\cdot\|_{p}$ is the $\ell_p$-norm on $\C^{n+1}$.
We shall give a summary of the various estimators of the family in the end of this section. 
Optimization problems associated to all of them rest upon a common principle -- minimization of the residual~$\| F_n[y - \vphi * y]\|_{p}$, with~$p \in \{2,\infty\}$,
regularized via the $\ell_1$-norm of the DFT of the filter.

The statistical properties of adaptive convolution-type estimators have been extensively studied. 
In particular, such estimators were shown to be nearly minimax-optimal, with respect to the pointwise loss and $\ell_2$-loss, for signals belonging to arbitrary, and unknown, shift-invariant linear subspaces of $\C(\Z)$ with bounded dimension, or sufficiently close to such subspaces as measured by the local $\ell_p$-norms, see~\cite{nemirovsky1991nonparametric,jn1-2009,jn2-2010,harchaoui2015adaptive,ostrovsky2016structure}. 
We give a summary of statistical properties of convolution-type estimators in Appendix~\ref{sec:background-adapt-est}.

However, the question of the algorithmic implementation of such estimators remains largely unexplored; in fact, we are not aware of any publicly available implementation of these estimators. 
Our goal here is to close this gap.
Note that problems similar to~\eqref{opt:l2pen-0} belong to the general class of second-order cone problems, and hence can in principle be solved to high numerical accuracy in polynomial time via interior-point methods~\cite{nemirovski2001complexity}. However, the computational complexity of interior-point methods grows polynomially with the problem dimension, and becomes prohibitive in signal and image denoising problems (for example, in image denoising this number is proportional to the number of pixels which might be as large as $10^8$). Furthermore, it is unclear whether high-accuracy solutions are necessary when the optimization problem is solved with the goal of obtaining a statistical estimator. In such cases, the level of accuracy sought, or the amount of computations performed, should rather be \textit{adjusted} to the statistical performance of the exact estimator itself. 
While these matters have previously been investigated in the context of linear regression~\cite{pilanci2016iterative} and~sparse recovery~\cite{bruer2015designing}, our work studies them in the context of convolution-type estimators. 

Notably, \eqref{opt:l2pen-0} and its counterparts have favorable properties:
\begin{itemize}
\item[-]
\emph{Easily accessible first-order information.} 
The objective value and gradient at a given point can be computed in time~$\cO(n\log n)$ via a series of Fast Fourier Transforms (FFT) and elementwise vector operations. 
\item[-]
\emph{Simple geometry.} 
After a straightforward re-parametrization, one is left with $\ell_1$-norm penalty or $\ell_1$-ball as a feasible set in the constrained formulation.
Prox-mappings for such problems, with respect to both the Euclidean and the ``$\ell_1$-adapted'' distance-generating functions, can be computed efficiently.
\item[-]
\emph{Medium accuracy is sufficient.}
We show that approximate solutions with specified (medium) accuracy preserve the statistical performance of the exact solutions.
\end{itemize}
All these properties make first-order optimization algorithms the tools of choice to deal with~\eqref{opt:l2pen-0} and similar problems.

\paragraph{Outline.}
In Section~\ref{sec:algos-general}, we recall two general classes of optimization problems, \textit{composite minimization}~\cite{beck2009fast,ActaNumerica2013} and \textit{composite saddle-point} problems~\cite{optbook2,ActaNumerica2013}, and the first-order optimization algorithms suitable for their numerical solution.
In Section~\ref{sec:implementation}, we show how to recast the optimization problems related to convolution-type estimators in one of the above general forms.
We then describe how to compute first-order oracles in the resulting problems efficiently using FFT.
In Section~\ref{sec:algos-complexity}, we establish problem-specific worst-case complexity bounds for the proposed first-order algorithms. 
These bounds are expressed in terms of the quantities that control the statistical difficulty of the signal recovery problem: signal length~$n$, noise variance~$\sigma^2$, and parameter~$r$ corresponding to the~$\ell_1$-norm of the Discrete Fourier tranform of the optimal solution.
A remarkable consequence of these bounds is that just~$\tilde\cO(\PSNR + 1)$ iterations of a suitable first-order algorithm are sufficient to match the statistical properties of an exact estimator; here~$\PSNR := \|F_{2n}[x]_{-n}^n\|_\infty / \sigma$ is the peak signal-to-noise ratio in the Fourier domain. 
This gives a rigorous characterization (in the present context) of the performance of ``early stopping'' strategies that allow to stop an optimization algorithm much earlier than dictated purely by the optimization analysis.
In Section~\ref{sec:experiments}, we present numerical experiments on simulated data which complement our theoretical analysis\footnotemark.
\footnotetext{The code reproducing all our experiments is available online at \url{https://github.com/ostrodmit/AlgoRec}.}

\paragraph{Notation.}
We denote $\C(\Z)$ the space of all complex-valued signals on $\Z$, or, simply, the space of all two-sided complex sequences.
We call $\C_n(\Z)$ the finite-dimensional subspace of $\C(\Z)$ consisting of signals supported on $[0, n]$:
\[
\C_n(\Z) = \left\{(x_\tau) \in \C(\Z): x_\tau = 0 \; \text{ whenever} \;  \tau \notin [0,n] \right\};
\]
its counterpart $\C_n^{\pm}(\Z)$ consists of all signals supported on $[-n, n]$. 
The unknown signal is assumed to come from one of such subspaces, which corresponds to a finite signal length.
Note that signals from $\C(\Z)$ can be naturally mapped to column vectors by means of the index-restriction operator $[\cdot]_{m}^n$, defined for any $m, n \in \Z$ such that $m \le n$ as
\[
[x]_{m}^n \in \C^{n-m+1}.
\]
In particular, $[\cdot]_{0}^n$ and $[\cdot]_{-n}^n$ define one-to-one mappings $\C_n(\Z) \to \C^{n+1}$ and $\C_n^\pm(\Z) \to \C^{2n+1}$.
For convenience, column-vectors in $\C^{n+1}$ and $\C^{2n+1}$ will be indexed starting from zero. 
We define the scaled $\ell_p$-seminorms on $\C(\Z)$:
\[
\| x \|_{n,p} := \frac{\|[x]_{0}^n\|_p}{(n+1)^{1/p}} = \left(\frac{1}{n+1}\sum_{\tau = 0}^n |x_\tau|^{p} \right)^{{1}/{p}}, \;  p \ge 1.
\]
We use the ``Matlab notation'' for matrix concatenation: $[A; B]$ is the vertical, and $[A, B]$ the horizontal concatenation of two matrices with compatible dimensions.
We introduce the unitary Discrete Fourier Transform (DFT) operator $F_n$ on $\C^{n+1}$, defined by
\begin{align*}
\displaystyle
[F_{n} x]_k = \frac{1}{\sqrt{n+1}} \sum\limits_{t = 0}^n x_{t} \exp\left({\frac{2\pi i k t}{n+1}}\right), \quad 0 \le k \le n.
\end{align*}
The unitarity of $F_n$ implies that its inverse $F_n^{-1}$ coincides with its conjugate transpose $F_n^\H$.
Slightly abusing the notation, we will occasionally shorten $F_n [x]_0^n$ to $F_n[x]$. 
In other words, $F_n[\cdot]$ is a map $\C_n(\Z) \to \C^{n+1}$, and the adjoint map $F_n^\H[x]$ simply sends $F_n^\H [x]_0^n$ to $\C_n(\Z)$ via zero-padding. 
We use the ``Big-O'' notation: for two non-negative functions $f, g$ on the same domain, $g = \cO(f)$ means that there is a generic constant $C \ge 0$ such that $g \le Cf$ for any admissible value of the argument; $g = \tilde\cO(f)$ means that $C$ is replaced with $C(\log^{\kappa}(n)+1)$ for some $\kappa > 0$; hereinafter $\log(\cdot)$ is the natural logarithm, and $C$ is a generic constant.

\paragraph{Estimators.} 
We now summarize all the estimators that are of interest in this paper. For brevity, we use the notation 
\begin{equation}
\label{eq:intro-residual}
\Res_p(\vphi) := \| F_n[y - \vphi * y]\|_{p}.
\end{equation}
\begin{itemize}
\item
\textit{{Con}strained {U}niform-{F}it estimator,} given for $\overline{r} \ge 0$ by
\begin{equation}
\label{opt:l8con}
\wh\vphi \in \Argmin_{\vphi\in \Phi_n(\overline r)} \quad
\Res_{\infty}(\vphi),
\tag{Con-UF}
\end{equation}
\[
\Phi_n(\overline r) := \left\{\vphi \in \C_n(\Z): \|F_n[\vphi]\|_{1} \leq {\overline{r}\over\sqrt{n+1}} \right\};
\]
\item 
\textit{{Con}strained Least-Squares estimator:}
\begin{equation}
\label{opt:l2con}
\wh\vphi \in \Argmin_{\vphi\in \Phi_n(\overline r)} \quad
\half \Res_{2}^2(\vphi);
\tag{Con-LS}
\end{equation}
\item
\textit{{Pen}alized {U}niform-Fit estimator:}
\begin{equation}
\label{opt:l8pen}
\wh \vphi \in \Argmin\limits_{\vphi \in \C_n(\Z)}
\Res_{\infty}(\vphi) + \lambda\|F_n[\vphi]\|_{1}; 
\tag{Pen-UF}
\end{equation}
\item 
\textit{{Pen}alized Least-Squares estimator:} 
\begin{equation}
\label{opt:l2pen}
\wh \vphi \in \Argmin\limits_{\vphi \in \C_n(\Z)}
\half \Res_{2}^2(\vphi) + \lambda\|F_n[\vphi]\|_{1}.
\tag{Pen-LS}
\end{equation}
\end{itemize}
We also consider~(\ref*{opt:l2con}$^*$) and~(\ref*{opt:l2pen}$^*$) -- counterparts of \eqref{opt:l2con} and~\eqref{opt:l2pen} in which $\half \Res_2^2(\vphi)$ is replaced with non-squared residual $\Res_2(\vphi)$. 
Note that~(\ref*{opt:l2con}$^*$) is equivalent to~\eqref{opt:l2con}, \ie~results in the same estimator; however, this does \textit{not} hold for~(\ref*{opt:l2pen}$^*$) and~\eqref{opt:l2pen}.


\section{Tools from Convex Optimization} 

In this section, we recall the tools from first-order convex optimization to be used later.
We describe two general types of optimization problems, \textit{composite minimization} and \textit{composite saddle-point} problems, together with efficient first-order algorithms for their solution. 
Following~\cite{ActaNumerica2013}, we begin by introducing the concept of \textit{proximal setup} which underlies these algorithms.

\label{sec:algos-general}

\subsection{Proximal Setup}
\label{sec:prox-setup}

Let a \textit{domain}~$U$ be a closed convex set in a Euclidean space $E$. A \textit{proximal setup} for $U$ is given by a norm $\|\cdot\|$ on $E$ (not necessarily Euclidean), and a \textit{distance-generating function} (d.-g.~f.) $\omega(u): U \to \R$, such that $\omega(u)$ is continuous and convex on $U$, admits a continuous selection $\omega'(u) \in \partial\omega(u)$ of subgradients on the set $\{ u \in U : \partial\omega(u) \ne \emptyset \}$, and is $1$-strongly convex with respect to~$\|\cdot\|$. 

The concept of proximal setup gives rise to several notions (see~\cite{ActaNumerica2013} for a detailed exposition): the $\omega$-center $u_\omega$, the Bregman divergence $D_u(\cdot)$, the $\omega$-radius $\Omega[\cdot]$ and the prox-mapping $\prox_{u}(\cdot)$ defined as
\begin{equation*}
\prox_{u}(g) = 
\argmin_{\xi \in U} \left\{\langle g, \xi \rangle + D_u(\xi) \right\} \, .
\end{equation*}

\paragraph{Blockwise Proximal Setups.}
We now describe a specific family of proximal setups which proves to be useful for our purposes. 
Let $E = \R^{N}$ with $N = 2(n+1)$; note that we can identify this space with $\C^{n+1}$ via (Hermitian) vectorization map $\Vec_n: \C^{n+1} \to \R^{2(n+1)}$,
\begin{equation}
\label{eq:vectorization}
\Vec_n z = [\Re(z_0); \, \Im(z_0); \, ...; \, \Re(z_n); \, \Im(z_n)].
\end{equation}
Now, supposing that
$
N = k(m+1)
$
for some non-negative integers $m, k$, let us split $u =[u^0; ...; u^m] \in \R^N$ into
$m+1$ blocks of size $k$, 
and equip $\R^N$ with the group $\ell_1/\ell_2$-norm:
\begin{equation}
\label{eq:block-norm}
\|u\| := \sum_{j=0}^m \|u^j\|_2.
\end{equation}
We also define the balls $U_N(R) := \{u \in \R^N: \|u\| \le R\}$.
\begin{theorem}[\cite{ActaNumerica2013}]
\label{th:prox-dgf}
Given~$E = \R^{N}$ as above, $\omega: \R^N \to \R$ defined by
\begin{equation}
\label{eq:prox-unconstrained}
\omega(u) = \frac{(m+1)^{(\tilde q-1)(2- \tilde q)/\tilde q}}{2\tilde c} \left[\sum_{j = 0}^m \|u^j\|_2^{\tilde q} \right]^{{2}/{\tilde q}}
\end{equation}
\begin{equation*}
\text{with} \; (\tilde{q}, \tilde{c}) = \left\{
\begin{array}{ll}
\left(2, \frac{1}{m+1}\right), &m \le 1,\\
\left(1 + \frac{1}{\log(m+1)}, \frac{1}{e\log(m+1)}\right), &m \ge 2,
\end{array}
\right.
\end{equation*}
is a d.-g.~f. for any ball $U_N(R)$ of the norm~\eqref{eq:block-norm} with $\omega$-center $u_\omega = 0$.
Moreover, for some constant $C$ and any $R \ge 0$ and $m, k \in \Z_+$, $\omega$-radius of $U_N(R)$ is bounded as
\begin{equation}
\label{eq:prox-radius-bound}
\Omega[U_N(R)] \le C(\sqrt{\log(m+1)}+1)R.
\end{equation}
\end{theorem}

We will use two particular cases of the above construction.
\begin{enumerate}
\item[(i)]
Case $m = n$, $k = 2$ corresponds to the $\ell_1$-norm on $\C^{n+1}$, and specifies the \textit{complex $\ell_1$-setup}. 
\item[(ii)]
Case $m = 0$, $k = N$ corresponds to the $\ell_2$-norm on $\C^{n+1}$, and specifies the \textit{$\ell_2$-setup} $(\|\cdot\|_2, \half\|\cdot\|_2^2)$.
\end{enumerate}
To work with them, we introduce specific norms on $\R^N$:
\begin{equation}
\label{eq:alg:u-norm}
\| u \|_{\C,p} := \|\Vec_n^{-1} u\|_p = \|\Vec_n^{\H}u\|_p, \quad p \ge 1.
\end{equation}
Note that $\| \cdot \|_{\C,1}$ gives the norm $\|\cdot\|$ in the complex $\ell_1$-setup, while $\| \cdot \|_{\C,2}$ coincides with the standard $\ell_2$-norm on $\R^N$.

\subsection{Composite Minimization Problems} 
\label{sec:fgm-general}

The general \textit{composite minimization} problem has the form
\begin{equation}
\label{eq:composite-min}
\min_{u \in U} \left\{ \phi(u) = f(u) + \Psi(u) \right\}.
\end{equation}
Here, $U$ is a domain in $E$ equipped with $\|\cdot\|$, $f(u)$ is convex and continuously differentiable on $U$, and $\Psi(u)$ is convex, lower-semicontinuous, finite on the relative interior of $U$, and can be non-smooth. 
Assuming that $U$ is equipped with a proximal setup $(\|\cdot\|, \omega(\cdot))$, let us define the \textit{composite prox-mapping}, see~\cite{beck2009fast}, as follows:
\begin{equation}
\label{eq:composite-prox}
\prox_{\Psi, u}(g)
= \argmin_{\xi \in U} \left\{ \langle g, \xi \rangle + D_u(\xi) + \Psi(\xi) \right\}.
\end{equation}

\paragraph{Fast Gradient Method.}
Fast Gradient Method (FGM), summarized as Algorithm~\ref{alg:fgm}, was introduced in~\cite{nesterov2007gradient} as an extension of the celebrated Nesterov algorithm for smooth minimization~\cite{nesterov1983method} to the case of constrained problems with non-Euclidean proximal setups.
It is guaranteed to find an approximate solution of~\eqref{eq:composite-min} with $\cO(1/T^2)$ accuracy after $T$ iterations. We defer the rigorous statement of this accuracy bound to Sec.~\ref{sec:algos-complexity}.


\begin{algorithm}[t]
\caption{Fast Gradient Method}
\label{alg:fgm}
\begin{algorithmic}
{
\REQUIRE stepsize $\eta > 0$\\
\STATE $u^0 = u_\omega$ \\
\STATE $g^0 =0 \in E$\\
\FOR{$t = 0, 1, ...$}
\STATE $u_t = \prox_{\eta \Psi, u_\omega}\left(\eta g^t\right)$ \\
\STATE $\tau_t = \frac{2(t+2)}{(t+1)(t+4)}$ \\
\STATE $u_{t+\frac{1}{3}} = \tau_t u_t + (1-\tau_t) u^t$ \\
\STATE $g_{t} = \frac{t+2}{2} \nabla f(u_{t+\frac{1}{3}})$ \\
\STATE $u_{t+\frac{2}{3}} = \prox_{\eta \Psi, u_t}\left(\eta g_{t}\right)$ \vspace{0.05cm}\\
\STATE $u^{t+1} = \tau_t u_{t+\frac{2}{3}} + (1-\tau_t) u^t$\\
\STATE $g^{t+1} = \sum_{\tau=0}^t g_{t}$\\
\ENDFOR
}
\end{algorithmic}
\end{algorithm}

\subsection{Composite Saddle-Point Problems}
\label{sec:cmp-general}

We also consider general~\textit{composite saddle-point} problems:
\begin{equation}
\label{eq:sp-gen}
\inf\limits_{u \in U}\, \max\limits_{v \in V} \left[\phi(u,v) = f(u,v) + \Psi(u)\right].
\end{equation}
Here, $U \subset E_u$ and  $V \subset E_v$ are domains in the corresponding Euclidean spaces $E_u, E_v$, and in addition~$V$ is compact; function $f(u,v)$ is convex in $u$, concave in $v$, and differentiable on $W := U \times V$; function $\Psi(u)$ is convex, lower-semicontinuous, can be non-smooth, and is such that $\prox_{\Psi, u}(g)$ is easily computable.
We can associate with $f$ a smooth vector field $F: W \to E_u \times E_v$, given by
\[
\begin{aligned}
&F([u; v]) = [\nabla_u f(u,v);  -\nabla_v f(u,v)].
\end{aligned}
\]
Saddle-point problem~\eqref{eq:sp-gen} specifies two convex optimization problems: that of minimization of $\overline\phi(u) = \max_{v \in V} \phi(u,v)$, or the primal problem, and that of maximization of $ \underline\phi(v) = \inf_{u \in U} \phi(u,v)$, or the dual problem.
Under the general conditions which hold in the described setting, see \eg~\cite{sion1958}, \eqref{eq:sp-gen} possesses an optimal solution $w^* = [u^*; v^*]$, called a \textit{saddle point}, such that the value of~\eqref{eq:sp-gen} is
$\phi(u^*,v^*) = \overline\phi(u^*) = \underline\phi(v^*)$, and $u^*, v^*$ are optimal solutions to the primal and dual problems.
The quality of a candidate solution $w = [u; v]$ can be evaluated via the \emph{duality gap} -- the sum of the primal and dual accuracies:
\begin{align*}
\overline\phi(u) - \underline\phi(v) = [\overline\phi(u) - \overline\phi(u^*)] + [\underline\phi(v^*) - \underline\phi(v)].
\end{align*}

\paragraph{Constructing the Joint Setup.}
When having a saddle-point problem at hand, one usually begins with ``partial'' proximal setups $(\|\cdot\|_{U}, \omega_{U})$ for $U \subseteq E_u$, and $(\|\cdot\|_{V}, \omega_V)$ for $V \subset E_v$, and must construct a ``joint'' proximal setup on $W$. Let us introduce the segment $U_* = [u^*, u_\omega]$, where $u_\omega$ is the $u$-component of the~$\omega$-center $w_\omega$ of $W$. Moreover, folllowing~\cite{ActaNumerica2013}, let us assume that the dual $\omega$-radius $\Omega[V]$ and the ``effective'' primal $\omega$-radius, defined as
\[
\Omega_*[U] := \min(\Omega[U], \Omega[U_*]),
\] 
are known (note that~$\Omega[U]$ can be infinite but $\Omega_*[U]$ cannot).
We can then construct a proximal setup
\begin{equation}
\label{eq:cmp-joint}
\begin{aligned}
&\| w \|^2 = \Omega^2[V] \, \|u\|_{U}^2 + \Omega^2_*[U] \, \|v\|_{V}^2, \\
&\omega(w) = \Omega^2[V] \, \omega_U(u) + \Omega^2_*[U] \, \omega_V(v).
\end{aligned}
\end{equation}
Note that the corresponding joint prox-mapping is reduced to the prox-mappings for the primal and dual setups.

\paragraph{Composite Mirror Prox.}
Composite Mirror Prox (CMP), introduced in~\cite{ActaNumerica2013} and summarized here as Algorithm~\ref{alg:cmp}, solves
the general composite saddle-point problem~\eqref{eq:sp-gen}. 
When applied with proximal setup~\eqref{eq:cmp-joint}, this algorithm admits an $\cO(1/T)$ accuracy bound after $T$ iterations; 
the formal statement is deferred to Sec.~\ref{sec:algos-complexity}.

\begin{algorithm}[t]
\caption{Composite Mirror Prox}
\label{alg:cmp}
\begin{algorithmic}
\REQUIRE stepsize $\eta > 0$\\
\STATE $w_0 := [u_0; v_0] = w_\omega$
\FOR{$t = 0, 1, ...$}
\STATE $w_{t+\frac{1}{2}} = \text{Prox}_{\eta \Psi, w_t}(\eta F(w_t))$
\STATE $w_{t+1} = \text{Prox}_{\eta \Psi, w_t}(\eta F(w_{t+\frac{1}{2}}))$
\STATE $w^{t+1} := [u^{t+1}; v^{t+1}] = \frac{1}{t+1}\sum_{\tau = 0}^{t} w_\tau$ 
\ENDFOR
\end{algorithmic}
\end{algorithm}

\section{Algorithmic Implementation}
\label{sec:implementation}



\paragraph{Change of Variables.}
When working with convolution-type estimators, our first step is to transfer the problem to the Fourier domain, so that the feasible set and the penalization term become quasi-separable.
Namely, noting that the adjoint map of $\Vec_n: \C^{n+1} \to \R^{2n+2}$, cf.~\eqref{eq:vectorization}, is given by
\[
\begin{aligned}
\Vec_n^\H u = [u_0; \, u_2; \, ...; \, u_{2n}] + i [u_1; \, u_3; \, ...; \, u_{2n+1}],
\end{aligned}
\]
consider the transformation
\begin{equation}
\label{eq:change_of_vars}
u = \Vec_n F_n[\varphi] \quad
b = \Vec_n F_n[y]
\end{equation}
Note that $\varphi = F_n^\H[\Vec_n^\H u] \in \C_n(\Z)$, and hence 
\[
\|F_n[y - y * \varphi]\|^2_{2} 
= \|Au - b\|_2^2,
\]
where $A: \R^{2n+2} \to \R^{2n+2}$ is defined by
\begin{equation}
\label{eq:operator-direct}
Au = \Vec_n F_n \left[y * F_n^\H[\Vec_n^\H u] \right].
\end{equation}


We are about to see that all recovery procedures can indeed be cast into one of the ``canonical'' forms~\eqref{eq:composite-min},~\eqref{eq:sp-gen}. Moreover, the gradient computation is then reduced to evaluating the convolution-type operator $A$ and its adjoint $A^\H = A^\T$.

\paragraph{Problem Reformulation.}
After the change of variables~\eqref{eq:change_of_vars}, problems~\eqref{opt:l2con} and~\eqref{opt:l2pen} 
take form~\eqref{eq:composite-min}:
\begin{align}
&\min\limits_{\|u\|_{\C,1} \le R}
\left[f(u) := \half\|Au - b\|_2^2\right]
+ \lambda \|u\|_{\C,1},
\label{eq:optim-fgd}
\end{align}
where $\|\cdot\|_{\C,p}$ is defined in~\eqref{eq:alg:u-norm}.
In particular, \eqref{opt:l2con} is obtained from~\eqref{eq:optim-fgd} by setting $\lambda = 0$ and~$R = \frac{\overline r}{\sqrt{n+1}},$
and~\eqref{opt:l2pen} is obtained by setting $R = \infty$.
Note that
\[
\nabla f(u) = A^\T (Au-b).
\]

On the other hand,  problems~\eqref{opt:l8con},~\eqref{opt:l8pen}, and~(\ref*{opt:l2con}$^*$)  
can be recast as saddle-point problems~\eqref{eq:sp-gen}. 
Indeed, the dual norm to $\|\cdot\|_{\C,p}$ is $\|\cdot\|_{\C,q}$ with~$q = \frac{p}{p-1}$, whence
\[
\|F_n[y - y * \varphi]\|_{p} = \|Au - b\|_{\C,p} = \max_{\|v\|_{\C,q} \le 1} \langle v, Au-b \rangle;
\]
as such, \eqref{opt:l8con}, \eqref{opt:l8pen} and~(\ref*{opt:l2con}$^*$) are reduced to a saddle-point problem
\begin{equation}
\label{eq:optim-mp}
\min\limits_{\|u\|_{\C,1} \le R} \;\; \max_{\|v\|_{\C,q} \le 1} 
[f(u,v) := \langle v, Au-b \rangle] 
+\lambda \|u\|_{\C,1},
\end{equation}
where~$q=1$ for~\eqref{opt:l8con} and~\eqref{opt:l8pen}, and~$q=2$ in case of~(\ref*{opt:l2con}$^*$). 
Note that $f(u,v)$ is bilinear, and one has
\[
[\nabla_u f(u,v); \nabla_v f(u,v)] = [A^\T v; Au-b].
\]
We are now in the position to apply the algorithms described in Sec.~\ref{sec:algos-general}.
One iteration of either of them is reduced to a few computations of the gradient (which, in turn, is reduced to evaluating $A$ and $A^\T$) and prox-mappings.
We now show how to evaluate operators $A$ and $A^\T$ in time~$\cO(n \log n)$.


\paragraph{Evaluation of $Au$ and $A^{\normalfont{\T}}v$.} 
Operator $A$, cf.~\eqref{eq:operator-direct}, can be evaluated in time $\cO(n\log n)$ via FFT.
The key fact is that the convolution $[y * \vphi]_0^n$ is contained in the first $n+1$ coordinates of the \textit{circular} convolution of $[y]_{-n}^n$ with a zero-padded filter 
$
\psi = [[\vphi]_0^n; 0_{n}] \in \C^{2n+1}.
$
Using the DFT diagonalization property, this fact can be expressed as
\begin{align*}
[y * \varphi]_t 
&= {\sqrt{2n+1}} \, [F^\H_{2n} D_y F_{2n} \psi]_t, \quad 0 \le t \le n,
\end{align*}
where operator $D_y = \textup{diag}(F_{2n} [y]_{-n}^n)$ on $\C^{2n+1}$ can be constructed in $\cO(n \log n)$ by FFT, and evaluated in $\cO(n)$.
Let $P_{n}: \C^{2n+1} \to \C^{n+1}$ project to the first $n+1$ coordinates of $\C^{2n+1}$; its adjoint $P_n^\H$ is the zero-padding operator which complements $[\vphi]_0^n$ with $n$ trailing zeroes. Then,
\begin{equation}
\label{eq:operator-long-exp}
Au = \sqrt{2n+1} \cdot \Vec_n^{\vphantom \H} F^{\vphantom \H}_n P^{\vphantom \H}_n F^\H_{2n} D_y F^{\vphantom \H}_{2n} P_n^\H F_n^\H \Vec_n^\H u,
\end{equation}
where all operators in the right-hand side can be evaluated in $\cO(n \log n)$.
Operator $A^\T = A^\H$ can be treated in the same manner by taking the adjoint of~\eqref{eq:operator-long-exp}.

\subsection{Computation of Prox-Mappings}
\label{sec:prox-computation}

It is worth mentioning that the composite prox-mappings in all cases of interest can be computed in time~$\cO(n)$; in some cases it can be done explicitly, and in others via a root-finding algorithm. 
These computations are described below.
It suffices to consider partial proximal setups separately; the case of joint setup in saddle-point problems can be treated using that the joint prox-mapping is separable in $u$ and $v$, cf.~Sec.~\ref{sec:cmp-general}.
Recall that the possible partial setups $(\|\cdot\|,\omega(\cdot))$ comprise the $\ell_2$-setup with $\|\cdot\| = \|\cdot\|_{\C,2} = \|\cdot\|_{2}$ and the (complex) $\ell_1$-setup
with $\|\cdot\| = \|\cdot\|_{\C,1}$; in both cases, $\omega(\cdot)$ is given by~\eqref{eq:prox-unconstrained}.
Computing $\prox_{\frac{1}{L}\Psi, u}(g)$, cf.~\eqref{eq:composite-prox}, amounts to solving 
\begin{equation}
\label{eq:prox-con}
\min_{\xi \in \R^{N}} \left\{ \xi^\T (g-\omega'(u)) + \omega(\xi): \|\xi\|_{\C,q} \le R\right\},
\end{equation}
in the constrained case, and
\begin{equation}
\label{eq:prox-pen}
\min_{\xi \in \R^{N}} \left\{ \xi^\T (g-\omega'(u)) + \omega(\xi) + \frac{\lambda}{L} \|\xi\|_{\C,1}^q \right\}, 
\end{equation}
in the penalized case\footnotemark;
\footnotetext{For the purpose of future reference, we also consider the case of squared $\|\cdot\|_{\C,1}$-norm penalty.}
in both cases, $q \in \{1,2\}$. In the constrained case with $\ell_2$-setup, the task is reduced to the Euclidean projection onto the $\ell_2$-ball if $q=2$, and onto the $\ell_1$-ball if $q=1$; the latter can be done (exactly) in $\tilde \cO(N)$ via the algorithm from~\cite{duchi2008efficient} -- for that, one first solves~\eqref{eq:prox-con} for the complex phases corresponding to the pairs of components of $\xi$.
The constrained case with $\ell_1$-setup is reduced to the penalized case by passing to the Langrangian dual problem. Evaluation of the dual function amounts to solving a problem equivalent to~\eqref{eq:prox-pen} with $q = 1$, and \eqref{eq:prox-con} can be solved by a simple root-finding procedure if one is able to solve~\eqref{eq:prox-pen}.
As for~\eqref{eq:prox-pen}, below we show how to solve it explicitly when $q = 1$, and reduce it to one-dimensional root search (so that it can be solved in $\cO(n)$ to numerical tolerance) when $q = 2$. 
%
%
%
%
Indeed,~\eqref{eq:prox-pen} can be recast in terms of  the complex variable $\zeta = \Vec_n^\H \xi$:
\begin{equation}
\label{opt:complex-prox}
\min_{\zeta \in \C^{n+1}} \left\{  \langle \zeta, z \rangle + \underline{\omega}(\zeta) + \frac{\lambda}{L} \|\zeta\|_1^q \right\}, 
\end{equation}
where 
$
z = \Vec_n^\H(g-\omega'(u)),
$
and $\underline{\omega}(\zeta) = \omega(\xi)$, cf.~\eqref{eq:prox-unconstrained}, whence
\begin{equation}\label{eq:adapted-prox-complex}
\underline{\omega}(\zeta) = \frac{C(m,\tilde q,\tilde \gamma)\|\zeta\|_{\tilde q}^2}{2},
\end{equation}
with $C(m,\tilde q,\tilde \gamma) = \frac{1}{\tilde \gamma} (m+1)^{(\tilde q-1)(2-\tilde q)/\tilde q}$.
Now,~\eqref{opt:complex-prox} can be minimized first with respect to the complex arguments, and then to the absolute values of the components of $\zeta$.
Denoting $\zeta^*$ a (unique) optimal solution of~\eqref{opt:complex-prox}, the first minimization results in
$
\zeta_j^* = -\frac{z_j}{|z_j|} |\zeta^*_j|, \; 0 \le j \le n,
$
and it remains to compute the absolute values $|\zeta^*_j|.$


\paragraph{Case $\boldsymbol{q=1}$.}
The first-order optimality condition implies
\begin{equation}
\label{eq:prox-optim-condition-q-1}
C(m,\tilde q,\tilde \gamma) \|\zeta^*\|_{\tilde q}^{2-\tilde q} |\zeta_j^*|^{\tilde q-1} + \frac{\lambda}{L} \mathds{1}\{ |\zeta_j^*| > 0 \} = |z_j|.
\end{equation}
Denoting $\tilde p = \frac{\tilde q}{\tilde q-1}$, and using the soft-thresholding operator 
\[
\Soft_M(x) = (|x|-M)_+\, \sign(x),
\] 
we obtain the explicit solution:
\begin{equation*}
\begin{aligned}
&\zeta_j^* = 
\frac{1}{C(m,\tilde q,\tilde \gamma)} \left( \frac{\theta_j}{\|\theta\|_{\tilde p}^{2-\tilde q}} \right)^{{\tilde p}/{\tilde q}}, \;\; \theta_j = \Soft_{\lambda / L}(z_j).
\end{aligned}
\end{equation*}
In the case of $\ell_2$-setup this reduces to $\zeta_j^* = \Soft_{\lambda / L}(z_j)$.

\paragraph{Case $\boldsymbol{q=2}$.}
Instead of~\eqref{eq:prox-optim-condition-q-1}, we arrive at
\begin{equation}
\label{eq:prox-optim-condition-q-2}
C(m,\tilde q,\tilde \gamma) \|\zeta^*\|_{\tilde q}^{2-\tilde q} |\zeta_j^*|^{\tilde q-1} + \frac{2\lambda \|\zeta^*\|_1}{L}  \mathds{1}\{ |\zeta_j^*| > 0 \} = |z_j|,
\end{equation}
which we cannot solve explicitly.
However, note that a counterpart of~\eqref{eq:prox-optim-condition-q-2}, in which $\|\zeta^*\|_1$ is replaced with parameter $t \ge 0$, can be solved explicitly similarly to~\eqref{eq:prox-optim-condition-q-1}. 
Let $\zeta^*(t)$ denote the corresponding solution for a fixed $t$, which can be obtained in $\cO(n)$ time.
Clearly, $\|\zeta^*(t)\|_1$ is a non-decreasing function on $\R_+$.
Hence, \eqref{eq:prox-optim-condition-q-2}  can be solved, up to numerical tolerance, by any one-dimensional root search procedure, in $\cO(1)$ evaluations of $\zeta^*(t)$.
\section{Theoretical Analysis}
\label{sec:algos-complexity}

The proofs of the technical statements of this section are collected in Appendix~\ref{sec:algos-proofs}.
\subsection{Bounds on Absolute Accuracy}
\label{sec:generic-acc}
We first recall from~\cite{ActaNumerica2013} the worst-case bounds on the \textit{absolute accuracy} in objective, defined as
$
\veps(t) := \phi(u^t) - \phi(u^*)
$
for composite minimization problems, and $\overline\veps(t) := \overline\phi(u^t) - \overline\phi(u^*)$ for saddle-point problems.
These bounds, summarized in~Theorems~\ref{th:fgm_perf}--\ref{th:cmp_perf} below, are applicable when solving \textit{arbitrary} problems of the types~\eqref{eq:composite-min}, \eqref{eq:sp-gen} with the suitable first-order algorithm, and are expressed in terms of the ``optimization'' parameters that specify the regularity of the objective and the $\omega$-radius.



\begin{theorem}
\label{th:fgm_perf}
Suppose that $f$ has $L_f$-Lipschitz gradient:
\[
\| \nabla f(u) - \nabla f(u') \|_* \le L_f \|u-u'\| \quad \forall u, u' \in U
\]
where $\|\cdot\|_*$ is the dual norm to $\|\cdot\|$,
and let $u^T$ be generated by $T$ iterations of Algorithm~\ref{alg:fgm} with stepsize $\eta = \frac{1}{L_f}$.
Then,
\[
\veps(T) = \cO\left(\frac{L_f \Omega_*^2[U]}{T^2}\right).
\]
\end{theorem}

\begin{theorem}
\label{th:cmp_perf}
Let $f(u,v)$ be as in~\eqref{eq:optim-mp}\footnote{For simplicity, we only state the bound for bilinear $f(u,v)$.},
and assume that vector field $F$ is $L_F$-Lipschitz on $W = U \times V$:
\[
\|F(w)-F(w')\|_* \le L_F \|w-w'\| \quad \forall w,w' \in W.
\]
Let $w^T = [u^T; v^T]$ be generated by $T$ iterations of Algorithm~\ref{alg:cmp} with joint setup~\eqref{eq:cmp-joint} and~$\eta = \frac{\Omega[V]}{\Omega_*[U]L_F}$.
Then,
\begin{align*}
\overline\veps(T) = \cO\left( \frac{L_F\Omega_*[U] \Omega[V]}{T} \right).
\end{align*}
\end{theorem}



Our next goal is to translate these bounds into the language of  ``statistical'' parameters such as the norm of exact estimator and the peak signal-to-noise ratio in the Fourier domain, cf. Sec.~\ref{sec:intro}.
Let us make a couple of observations beforehand.

The first observation concerns the proximal setups to be used, and allows to control the $\omega$-radii.
If the partial domain (for $u$ or $v$) is an $\|\cdot\|_{\C,2}$-norm ball, we will naturally use the $\ell_2$-setup in that variable. 
If the domain is an $\|\cdot\|_{\C,1}$-norm ball, we will consider choosing between the $\ell_1$-setup which is ``adapted'' to the geometry of the problem, see~\cite{ActaNumerica2013}, or the $\ell_2$-setup due to its simplicity in use.
Note that in all these cases, the partial domains either coincide with or are contained in the balls~$U_N(1), U_N(R)$ of the corresponding norms, cf.~\eqref{eq:prox-radius-bound}, whence $\omega$-radii $\Omega[V], \Omega_*[U]$ can be bounded as follows:
\begin{equation}
\label{app:fgm-rad}
\Omega[V] = \tilde\cO(1), \quad \Omega_*[U] = \tilde\cO\left({r}/{\sqrt{n+1}}\right),
\end{equation}
where 
\begin{equation}
\label{eq:opt-norm-parameter}
r = \sqrt{n+1} \|F_n[\wh\vphi]\|_1
\end{equation}
is the scaled norm of an optimal solution (note that $\overline{r} \ge r$).

The second observation concerns the Lipschitz constants $L_f, L_F$ in the chosen setups. It is convenient to define parameters $q_u, q_v$ that take values in $\{2,1\}$
depending on the partial setup used in the corresponding variable;
besides, let $p_u = \frac{q_u}{q_u-1}$ and~$p_v = \frac{q_v}{q_v-1}$.
Introducing the complex counterpart of $A$, operator $\A: \C^{n+1} \to \C^{n+1}$ given by
\[
\A [\vphi]_0^n = F_n [y* F_n^\H[\vphi]_0^n] \;\; \Leftrightarrow \;\; A = \Vec_n^{\vphantom\H} \circ \A \circ \Vec_n^\H,
\]
we can conveniently express Lipshitz constants~$L_f$, $L_F$ in terms of operator norms~$\|\A\|_{\alpha \to \beta} := \sup_{\|\psi\|_{\alpha} = 1} \|\A \psi \|_\beta$:
\begin{equation}
\label{eq:lip-bounds}
\begin{aligned}
\|\A\|_{1 \to 2}^2 \le L_f &= \|\A\|_{q_u \to 2}^2 \le \|\A\|_{2 \to 2}^2,\\
\|\A\|_{1 \to \infty} \le L_F &= \|\A\|_{q^{\vphantom{*}}_u \to p_v} \le \|\A\|_{2 \to 2}.
\end{aligned}
\end{equation}
Now, the norm~$\|\A\|_{2 \to 2}$ itself can be bounded as follows:
\begin{lemma}
\label{th:op-norm-upper}
One has
\[
\|\A\|_{2 \to 2} \le \sqrt{2n+1} \cdot \|F_{2n} [y]_{-n}^n\|_{\infty}.
\]
\end{lemma}
Together with~\eqref{app:fgm-rad}, Lemma~\ref{th:op-norm-upper} results in
\begin{prop}
\label{th:abs-acc}
Solving~\eqref{opt:l2con} or \eqref{opt:l2pen} by Algorithm~\ref{alg:fgm} with proximal setup as described above, one has
\begin{equation}
\label{eq:fgm-abs-acc}
\veps(T) = 
\tilde\cO\left(\frac{r^2\|F_{2n} [y]_{-n}^n\|_{\infty}^2}{T^2}\right).
\end{equation}
Similarly, solving~\eqref{opt:l8con}, \eqref{opt:l8pen}, \textup{(\ref*{opt:l2con}$^*$)}, or~\textup{(\ref*{opt:l2pen}$^*$)} by Algorithm~\ref{alg:cmp} with proximal setup as described above, 
\begin{equation}
\label{eq:cmp-abs-acc}
\overline\veps(T) = 
\tilde\cO\left(\frac{r\|F_{2n} [y]_{-n}^n\|_{\infty}}{T}\right).
\end{equation}
\end{prop}

\paragraph{Discussion: comparison of setups.}
Note that 
Proposition~\ref{th:abs-acc} gives the same upper bound on the accuracy $\veps(T)$ irrespectively of the chosen proximal setup. 
This is because we used the operator norm $\|\cA\|_{2 \to2}$ as an upper bound for~$L_f$ and~$\sqrt{L_F}$ while these quantities are in fact equal to~$\|\cA\|_{1 \to 2}$ or~$\|\cA\|_{1 \to \infty} \le \|\cA\|_{1 \to 2}$ when one uses the ``geometry-adapted'' $\ell_1$-setup in at least one of the variables. 
For a general linear operator $\cA$ on $\C^{n+1}$ the gaps between $\|\cA\|_{2 \to 2}$ and the latter norms can be as large as $\sqrt{n+1}$ or $n+1$, hence one might expect the bound of Proposition~\ref{th:abs-acc} to be loose.
However, intuitively~$\cA$ is ``almost'' a diagonal operator -- it would as such is we worked with the \textit{circular} convolution. Hence, we can expect its various $\|\cdot\|_{q \to p}$ norms in~\eqref{eq:lip-bounds} to be mutually close (in the case $\cA = \text{diag}(a)$ they all coincide with $\|a\|_\infty$).
This heuristic observation can be made precise:
\begin{prop}
\label{th:op-norm-lower}
Assume that $\sigma = 0$, and $x \in \C(\Z)$ is $(n+1)$-periodic:
$
x_\tau = x_{\tau-n-1}, \tau \in \Z.
$ 
Then, one has 
\[
\|\A\|_{1 \to \infty} = \sqrt{n+1} \|F_n[x]\|_{\infty}.
\]
\end{prop}

\subsection{Statistical Accuracy and  Complexity Bounds}
\label{sec:algos-stats}
In this section, we first characterize the \textit{statistical accuracy} of adaptive recovery procedures, defined as the absolute accuracy $\veps_*$ sufficient for the corresponding approximate estimator $\tilde \vphi$ to admit the same, up to a constant factor, theoretical risk bound as the exact estimator~$\wh \vphi$. The exact meaning of ``risk bound'' here depends on the estimator in consideration: for uniform-fit estimators it is the bound on the pointwise loss that was proved in~\cite{harchaoui2015adaptive}, and for least-squares estimators it is the bound on the~$\ell_2$-loss  proved in~\cite{ostrovsky2016structure}. 
The next two results state that statistical accuracy, defined in this sense, can be chosen as~$\sigma r$ for uniform-fit procedures, and~$\sigma^2r^2$ for least-squares procedures. 
The arguments, provided in Appendix~\ref{sec:algos-proofs}, closely follow those in~\cite{harchaoui2015adaptive} and~\cite{ostrovsky2016structure}.
\begin{theorem}
\label{th:stat-acc-l8}
An~$\veps_*$-accurate solution~$\tilde \vphi$ to~\eqref{opt:l8con} with~$\overline r = r$, or to~\eqref{opt:l8pen} with
$
\lambda = 16\sigma\sqrt{\left(n+1\right)\left(1+\log\left(\frac{n+1}{\delta}\right)\right)},
$ 
in both cases with $\veps_* = \cO(\sigma r)$, with prob.~$\ge 1-\delta$ satisfies
\begin{equation}
\label{eq:stat-acc-l8}
|x_n - [\tilde \vphi * y]_n| \le \frac{C \sigma r^2\sqrt{1+\log\left(\frac{n+1}{\delta}\right)}}{\sqrt{n+1}}.
\end{equation}
\end{theorem}
While Theorem~\ref{th:stat-acc-l8} controls the pointwise loss for uniform-fit estimators, the next theorem controls the $\ell_2$-loss for least-squares estimators.
To state it, we recall that a linear subspace $\S$ of $\C(\Z)$ is called \textit{shift-invariant} if it is an invariant subspace of the lag operator $\Delta$: $[\Delta x]_\tau = x_{\tau - 1}$ on $\C(\Z)$.
\begin{theorem}
\label{th:stat-acc-l2}
Assume that $x$ belongs to a shift-invariant subspace $\S$ with $\dim(\S) \le n$. 
Then, an $\veps_*$-accurate solution $\tilde \vphi$ to~\eqref{opt:l2con} with $\overline r = r$ or to~\eqref{opt:l2pen} with
\[
\lambda = 8\sqrt{2}\sigma^2\sqrt{n+1}\left(2+\log\left(\frac{8(n+1)}{\delta}\right)\right),
\]
in all cases with $\veps_* = \cO(\sigma^2 r^2)$, with prob.~$\ge 1-\delta$ satisfies
\begin{equation}
\label{eq:stat-acc-l2}
\|x - \tilde \vphi * y\|_{n,2} \le \frac{C\sigma \left(r\sqrt{1+\log \left(\frac{n+1}{\delta}\right)} + \sqrt{\dim(\S)}\right)}{\sqrt{n+1}}.
\end{equation}
\end{theorem}




\paragraph{Complexity Bound.}
Combining Theorems~\ref{th:stat-acc-l8}--\ref{th:stat-acc-l2} with Proposition~\ref{th:abs-acc}, we arrive at the following conclusion: for both classes of estimators, the number of iterations $T_*$ of the suitable first-order algorithm (Algorithm~\ref{alg:fgm} for the least-squares estimators and Algorithm~\ref{alg:cmp} for the uniform-fit ones) that guarantees accuracy $\veps_*$, with high probability satisfies
\begin{equation}
\label{eq:T-psnr}
T_* = \tilde\cO\left(\|F_{2n}[y]_{-n}^n\|_{\infty}/\sigma\right) = \tilde\cO\left(\PSNR + 1\right).
\end{equation}
Here, $\PSNR := \|F_{2n}[x]_{-n}^n\|_{\infty} / \sigma$ is the peak signal-to-noise ratio in the Fourier domain, and we used the unitary invariance of the complex Gaussian distribution.
Moreover, if it is known that the signal is sparse in the Fourier domain, that is, $\S$ is spanned by $s$ complex exponentials $e^{i\omega_k \tau}$ with frequencies on the grid, $\omega_k \in \left\{\frac{2\pi j}{n+1}, \; j \in \Z\right\}$, we can write 
\begin{equation}
\label{eq:psnr-sparse}
\PSNR = \cO(\SNR\sqrt{s})
\end{equation}
where $\SNR =  \|x\|_{n,2}/\sigma$ is the usual signal-to-noise ratio.

\begin{figure}[t]
\centering
\includegraphics[scale=0.5, clip=true, angle=0]{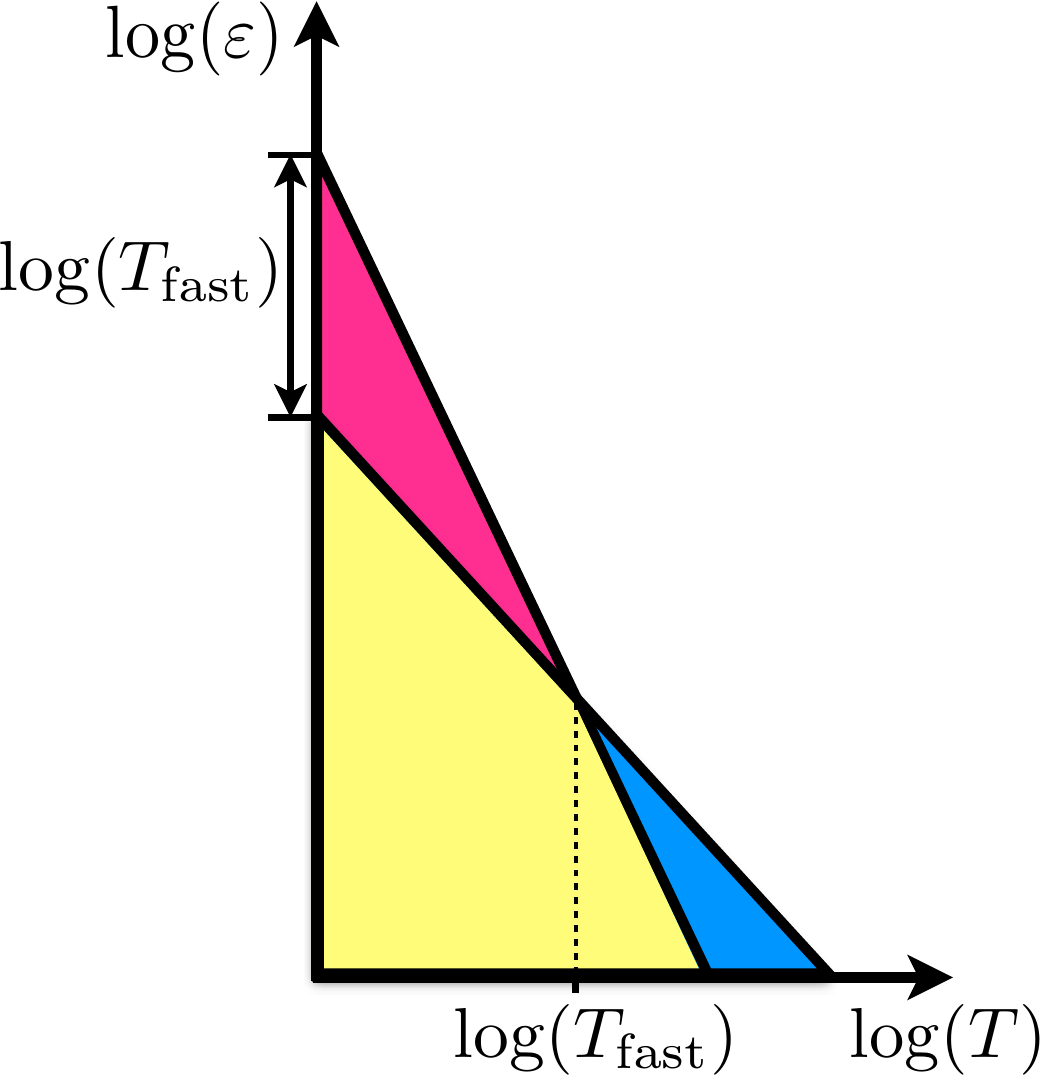}
\caption{``Phase transition'' for~Algorithm~\ref{alg:fgm}. The different slopes correspond to~\eqref{eq:fgm-root} and~\eqref{eq:fgm-slow}.}
\label{fig:elbow}
\end{figure}

\paragraph{Discussion: different ways of solving~\eqref{opt:l2con}.}
Note that Algorithm~\ref{alg:cmp} can be used to solve problems~(\ref*{opt:l2con}$^*$) and~(\ref*{opt:l2pen}$^*$) with non-squared residual by reducing them to (composite) saddle-point problems as shown in Sec.~\ref{sec:implementation}.
Hence, when solving~\eqref{opt:l2con} we have two alternatives: either to solve it directly with Algorithm~\ref{alg:fgm}, or to solve instead the equivalent problem~(\ref*{opt:l2con}$^*$) with Algorithm~\ref{alg:cmp}. 
Note that the complexity bound~\eqref{eq:T-psnr} only holds when Algorithm~\ref{alg:fgm}, and we can guess that this way of treating~\eqref{opt:l2con} is more beneficial.
Indeed, whenever the optimal residual~$\Res_2(\wh\vphi)$ is strictly positive, attaining accuracy $\sigma^2 r^2$ for~\eqref{opt:l2con} is equivalent to attaining accuracy
$
\veps_{**} = \frac{\sigma^2 r^2}{\Res_2(\wh\vphi)},
$ 
rather than $\veps_* = \sigma r$, for~(\ref*{opt:l2con}$^*$), where $\Res_2(\wh\vphi)$ is the optimal residual. Using Proposition~\ref{th:abs-acc}, the number of iterations of Algorithm~\ref{alg:cmp} to guarantee that is
\[
T_{**} = \frac{\Res_2(\wh\vphi)}{\sigma} \cdot \cO(\PSNR + 1).
\]
Potentially, this is much worse than~\eqref{eq:T-psnr} since $\Res_2(\wh\vphi)$ is expected to scale as the $\ell_2$-norm of the noise, \ie~$\sigma \sqrt{n+1}$. 


One curious property of Algorithm~\ref{alg:fgm} in the present context is its fast $\cO(1/T^{2})$ convergence in terms of the objective of~(\ref*{opt:l2con}$^*$). This fact, although surprizing at a first glance since the objective of~(\ref*{opt:l2con}$^*$) is \textit{non-smooth}, has a simple explanation. Note that in case of~\eqref{opt:l2con}, \eqref{eq:fgm-abs-acc} becomes
\begin{align}
\label{eq:fgm-squares}
\Res_2^2(\tilde\vphi) - \Res_2^2(\wh\vphi) = \tilde \cO\left(\frac{r^2 \|F_{2n} [y]_{-n}^n\|_{\infty}^2}{T^2}\right).
\end{align}
Dividing by~$\Res_2(\tilde\vphi) + \Res_2(\wh\vphi) \ge 2\Res_2(\wh\vphi)$, we obtain 
\begin{align}
\label{eq:fgm-root}
\Res_2(\tilde\vphi) - \Res_2(\wh\vphi) = \tilde \cO\left(\frac{r^2 \|F_{2n} [y]_{-n}^n\|_{\infty}^2}{\Res_2(\wh\vphi) T^2}\right),
\end{align}
\ie~$\cO(1/T^{2})$ convergence for~(\ref*{opt:l2con}$^*$) if~$\Res_2(\wh\vphi) > 0$. 
Moreover, this bound is crucial to achieve~\eqref{eq:T-psnr}, since~\eqref{eq:T-psnr} is exactly what is required for the right-hand side of~\eqref{eq:fgm-root} to be upper-bounded by $\veps_{**} = \frac{\sigma^2 r^2}{\Res_2(\wh\vphi)}$.

Finally, note that for small $T$, the $\cO(1/T^{2})$ bound~\eqref{eq:fgm-root} is dominated by the $\cO(1/T)$ bound
\begin{equation}
\label{eq:fgm-slow}
\Res_2(\tilde\vphi) - \Res_2(\wh\vphi) = \tilde \cO\left(\frac{r\|F_{2n} [y]_{-n}^n\|_{\infty}}{T}\right),
\end{equation}
which is obtained from~\eqref{eq:fgm-squares} by putting~$\Res_2(\wh\vphi)$ into the right-hand side and taking the square root. 
Hence, we can expect to see the faster $\cO(1/T^{2})$ convergence after 
\begin{align}
\label{eq:fgm-elbow}
T_{\fast} 
&= \frac{r \|F_{2n}[y]_{-n}^n\|_{\infty}}{\Res_2(\wh\vphi)}
\end{align}
iterations of Algorithm~\ref{alg:fgm}, as graphically shown in Fig.~\ref{fig:elbow}.
\section{Experiments} 
\label{sec:experiments}



\newcommand\halfwth{0.48\textwidth}
\newcommand\halfhgt{0.24\textheight}
\newcommand\thirdwth{0.32\textwidth}
\newcommand\thirdhgt{0.16\textheight}
\newcommand\quartwth{0.238\textwidth}
\newcommand\quarthgt{0.119\textheight}
\newcommand\sixthwth{0.156\textwidth}
\newcommand\sixthhgt{0.078\textheight}

\newcommand\compressedhgt{\quarthgt} 

\begin{figure}[]
\center
\begin{minipage}{\thirdwth}
\centering
\includegraphics[width=1\textwidth, height=\thirdhgt, clip=true, angle=0]{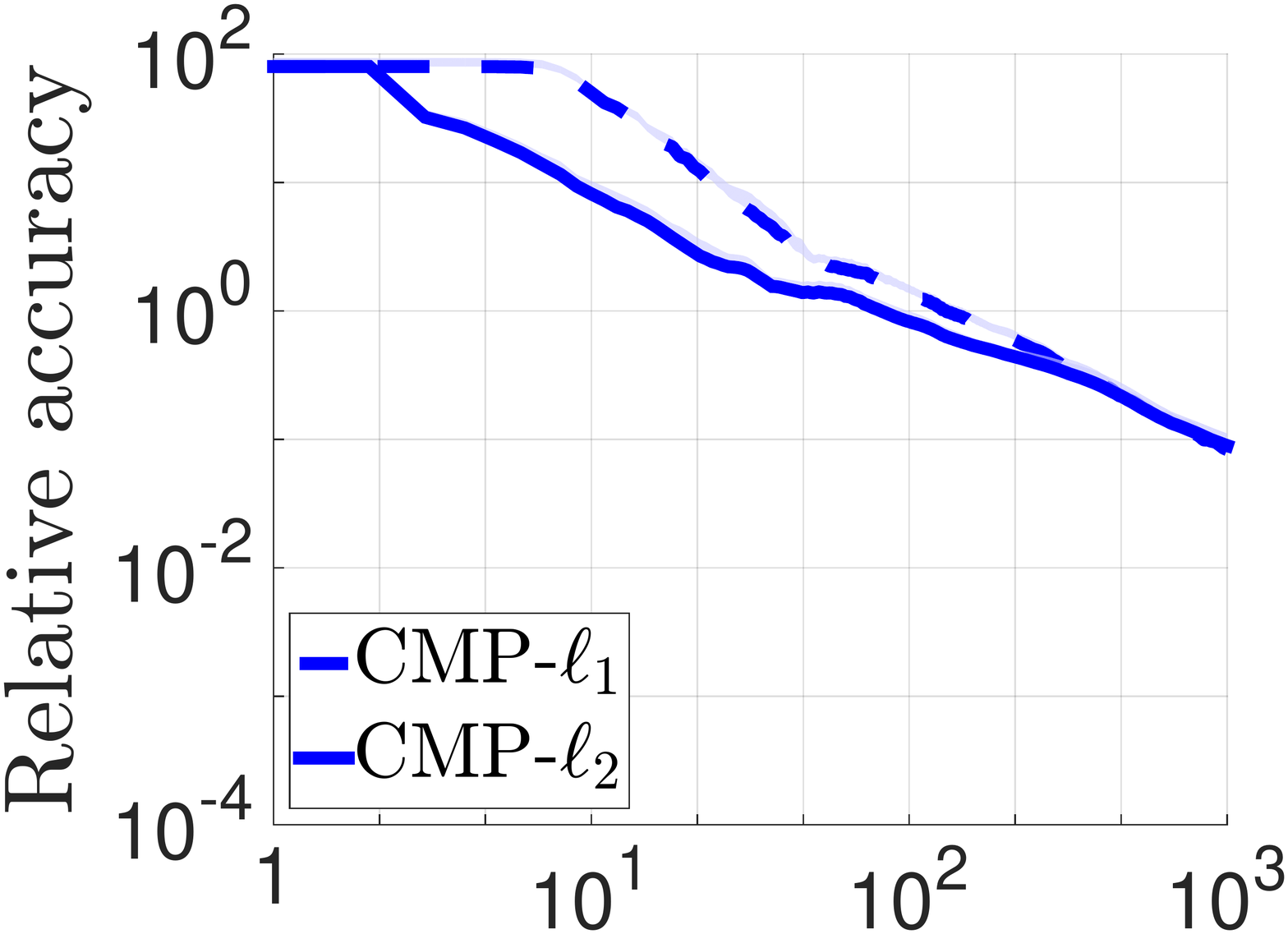}
\end{minipage}
\begin{minipage}{\thirdwth}
\centering
\includegraphics[width=1\textwidth, height=\thirdhgt, clip=true, angle=0]{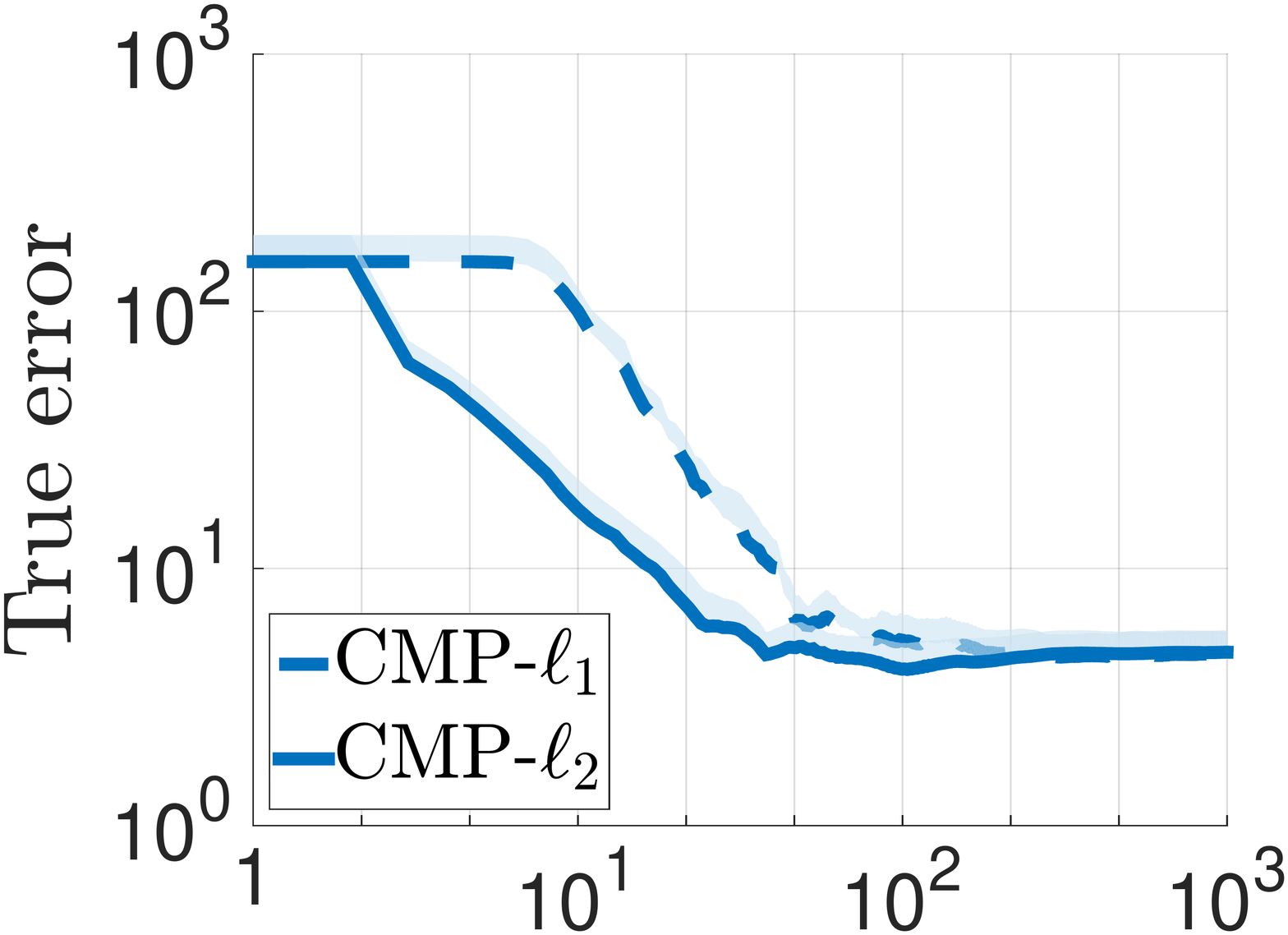}
\centering
\end{minipage}
\caption{Relative accuracy, left,
and~$\ell_\infty$-loss~$\|F_n[x -\tilde\vphi(T)* y]\|_{\infty}$, right, vs. iteration for approximate solutions to~\eqref{opt:l8con} by Algorithm~\ref{alg:cmp} in~\CohSin-8 with $\SNR = 16$.}
\label{fig:algos-exp1}
\end{figure}

In this series of experiments, our goal is to demonstrate the effectiveness of the approach and illustrate the theoretical results of Sec.~\ref{sec:algos-complexity}.
We estimate signals coming from an unknown shift-invariant subspace $\S$, implementing the following experimental protocol.
First, a random signal $[x_0; ...; x_{n}]$ with $n=100$ is generated according to one of the scenarios described below ($s$ is a parameter in both scenarios).
Then, $x$ is normalized so that $\|[x]_0^{n}\|_2 = 1$, and corrupted by i.i.d. Gaussian noise with a chosen level of $\SNR = (\sigma\sqrt{n})^{-1}$. 
A number of independent trials is performed to ensure the statistical significance of the results.
\begin{itemize}
\item In scenario \RanSin-$s$, the signal is a harmonic oscillation with $s$ frequencies: $x_t = \sum_{k=1}^s a_k e^{i\omega_k t}$.
The frequencies are sampled uniformly at random on $[0,2\pi[$, and the amplitudes uniformly on $[0,1]$.
\item In scenario \CohSin-$s$, we sample $s$ \textit{pairs} of close frequencies. Frequencies in each pair have the same amplitude and are separated only by $\frac{0.2\pi}{n}$ -- $0.1$ DFT bin -- so that the signal violates the usual frequency separation conditions, see~\eg~\cite{recht2}.
\end{itemize}
For constrained estimator we set $\overline r = 2\dim(\S)$ as suggested in~\cite{ostrovsky2016structure} for two-sided filters. 
Note that $\dim(\S) = s$ in \RanSin-$s$ and $\dim(\S) = 2s$ in \CohSin-$s$. 


\begin{figure}[t!]
\center
\begin{minipage}{\thirdwth}
\centering
\includegraphics[width=1\textwidth, height=\thirdhgt, clip=true, angle=0]{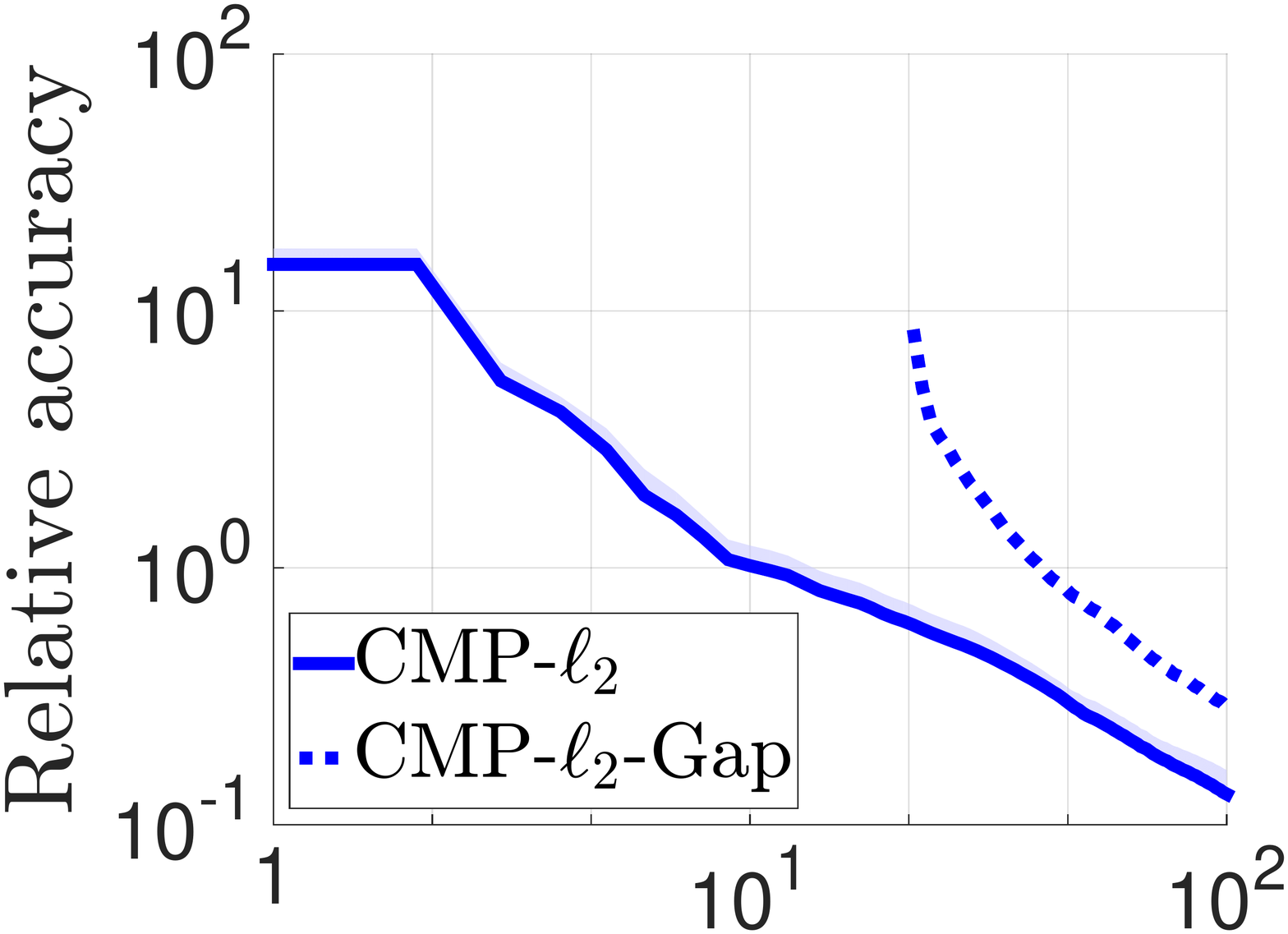}
\end{minipage}
\begin{minipage}{\thirdwth}
\centering
\includegraphics[width=1\textwidth, height=\thirdhgt, clip=true, angle=0]{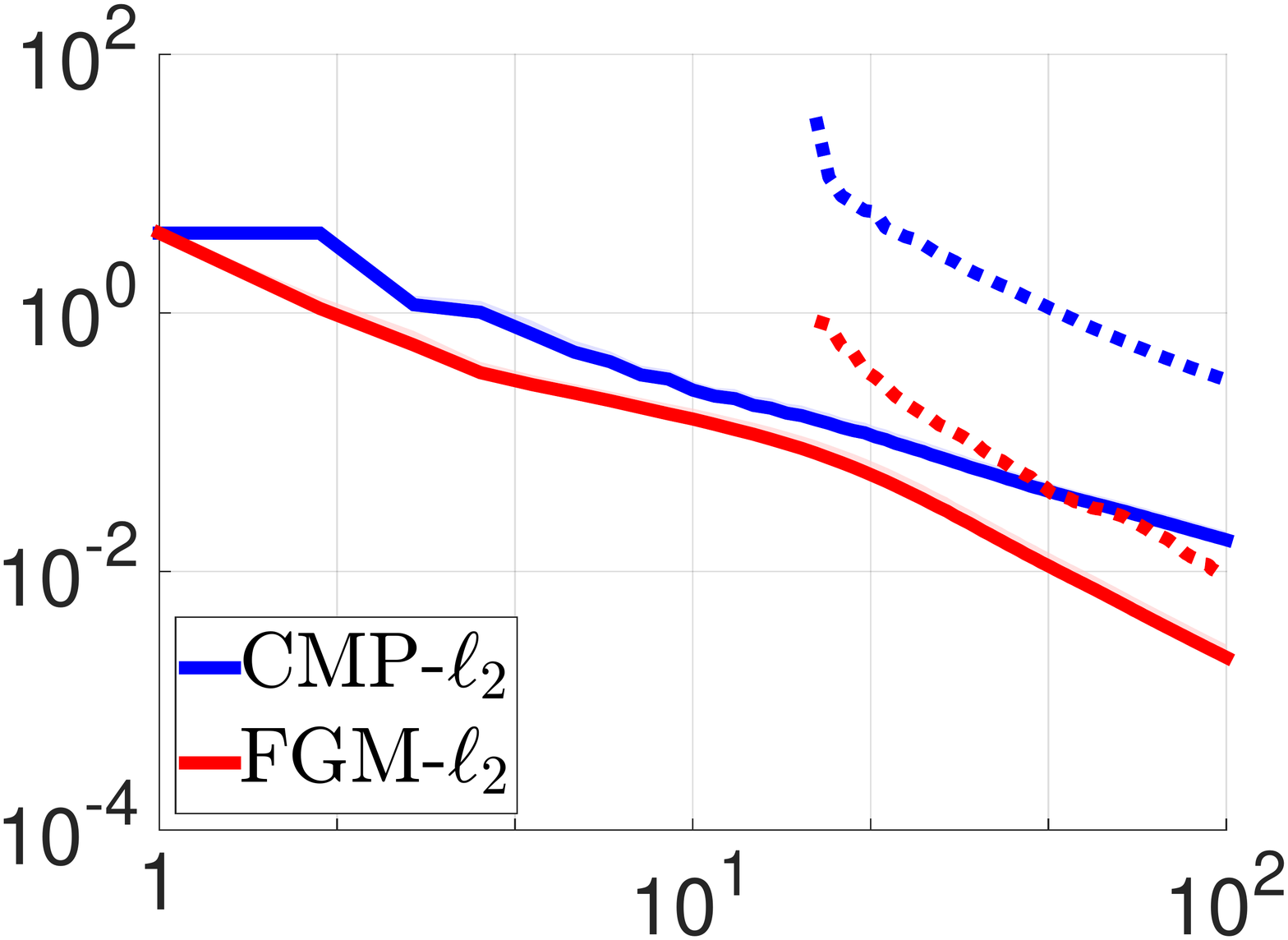}
\end{minipage}
\caption{Relative accuracy vs. iteration for~\eqref{opt:l8con}, left, and~(\ref*{opt:l2con}$^*$), right, in scenario \CohSin-4 with $\SNR = 4$. Dotted: accuracy certificates, see~\cite{nemirovski2010accuracy}.}

\label{fig:algos-exp3}
\end{figure}

\paragraph{Proof-of-Concept.}

In this experiment, we study estimator~\eqref{opt:l8con} in scenarios~\RanSin-$16$ and \CohSin-$8$. 
We run a version of CMP (Algorithm~\ref{alg:cmp}) with adaptive stepsize, see~\cite{ActaNumerica2013}, plotting the relative accuracy of the corresponding approximate solution $\tilde\vphi(T)$, that is, $\veps(T)$ normalized by the optimal value of the residual~$\Res_{\infty}(\wh\vphi)$, versus~$T$.
We also trace the true estimation error as measured by the $\ell_\infty$-loss in the Fourier domain,
$
\|F_n[x -\tilde\vphi(T)* y]\|_{\infty}.
$
Two joint proximal setups are considered: 
the full $\ell_2$-setup composed from the partial $\ell_2$-setups, and the full $\ell_1$-setup composed from the partial $\ell_1$-setups.
To obtain a proxy for~$\wh\vphi$, we recast~\eqref{opt:l8con} as a second-order cone problem, and run the MOSEK interior-point solver~ \cite{mosek}; note that this method is only available for small-sized problems.
We show upper $95\%$-confidence bounds for the convergence curves.

The results of this experiment, shown in Fig.~\ref{fig:algos-exp1}, can be summarized as follows. 
First, we see that the complexity of the optimization task grows with~$\SNR$ as predicted by~\eqref{eq:cmp-abs-acc}.
Second, provided that the number of frequencies is the same, there is no significant difference between scenarios \RanSin~and \CohSin~for the computational performance of our algorithms  (albeit we find $\CohSin$ to be slightly harder, and we only show the results for this scenario here).
We also find, somewhat unexpectedly, that the $\ell_2$-setup outperforms the ``geometry-adapted'' setup in earlier iterations; however, the performances of the two setups match in later iterations. 

Overall, we find that the first $100$ iterations result in~$100\%$ relative accuracy, \ie.~$\Res_{\infty}(\tilde\vphi) \le 2\Res_{\infty}(\wh\vphi)$.
In fact, from the analysis of uniform-fit estimators in the proof of Theorem~\ref{th:stat-acc-l8} we can derive the bound $\Res_{\infty}(\wh\vphi) = \tilde\cO(\sigma r)$, implying that the conditions of~Theorem~\ref{th:stat-acc-l8} are met for $\tilde\vphi$. 
As such, we can predict that further optimization is redundant. 
This is empirically confirmed: the true error begins to plateau after no more than~$100$ iterations.

\paragraph{Convergence and Accuracy Certificates.}

Here we illustrate the convergence of FGM (Algorithm~\ref{alg:fgm}) and CMP (Algorithm~\ref{alg:cmp}), including the case of~(\ref*{opt:l2con}$^*$) where both algorithms can be applied and thus compared.
We work in the same setting as previously, but this time also study~(\ref*{opt:l2con}$^*$) for which we compare the recommended approach via Algorithm~\ref{alg:fgm} and the alternative approach via~Algorithm~\ref{alg:cmp}  as discussed in Sec.~\ref{sec:algos-stats}.
The results are shown in Fig.~\ref{fig:algos-exp3}. 
We empirically observe~$\cO(1/T)$ convergence of Algorithm~\ref{alg:cmp} when solving~\eqref{opt:l8con}, as well as~$\cO(1/T^{2})$ convergence of Algorithm~\ref{alg:fgm} when solving~(\ref*{opt:l2con}$^*$), after a certain threshold as predicted by~\eqref{eq:fgm-root}--\eqref{eq:fgm-elbow}. 
In addition to accuracy curves, we plot upper bounds on them obtained via the technique of accuracy certificates, see~\cite{nemirovski2010accuracy} and Appendix~\ref{sec:algos-certificates}.
Such bounds can be used to stop the algorithms once the desired accuracy has been attained. 



\paragraph{Statistical Complexity Bound.}

\begin{figure}[t]
\center
\begin{minipage}{\thirdwth}
\centering
\includegraphics[width=1\textwidth, height=\thirdhgt, clip=true, angle=0]{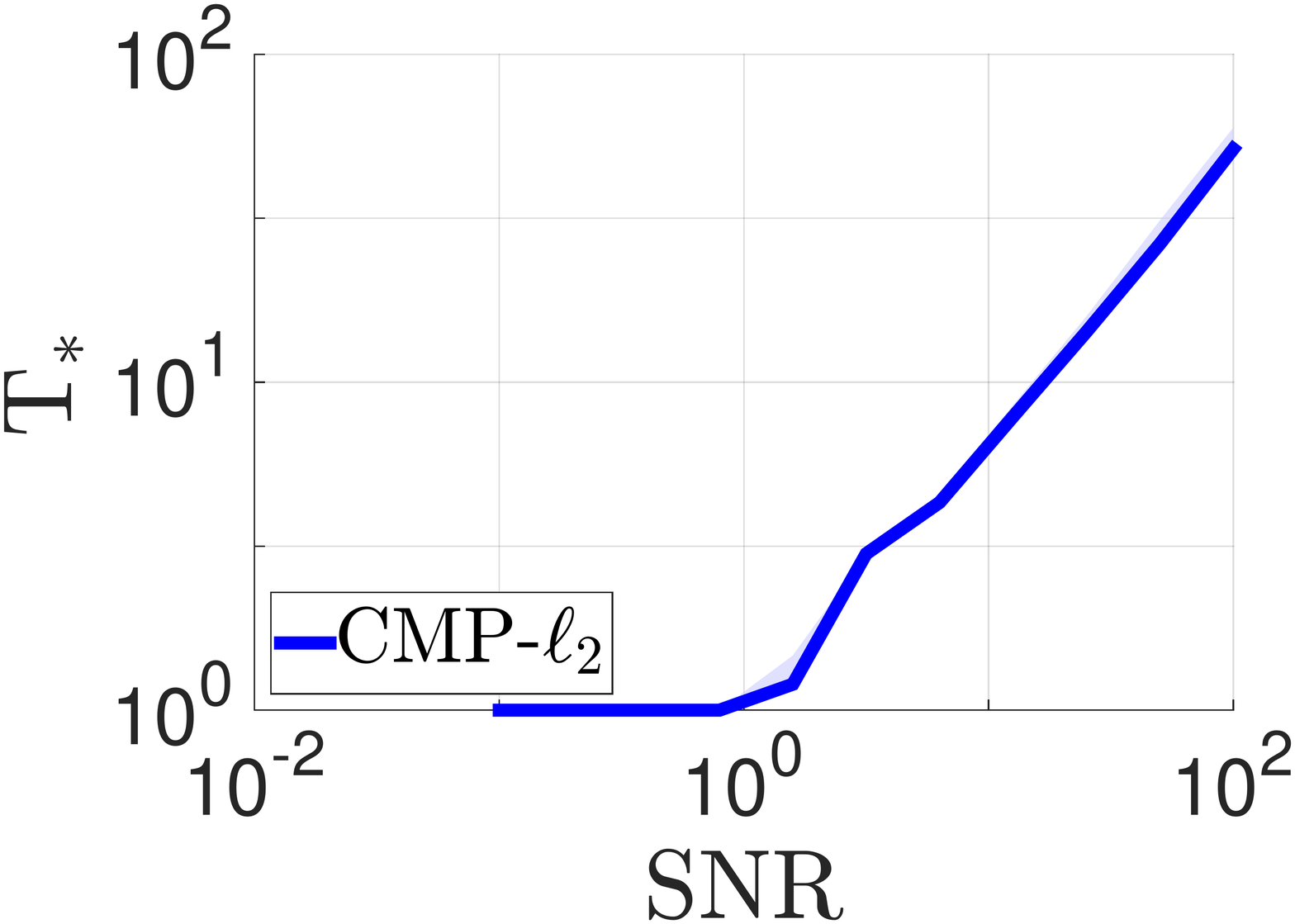}
\end{minipage}
\begin{minipage}{\thirdwth}
\centering
\includegraphics[width=1\textwidth, height=\thirdhgt, clip=true, angle=0]{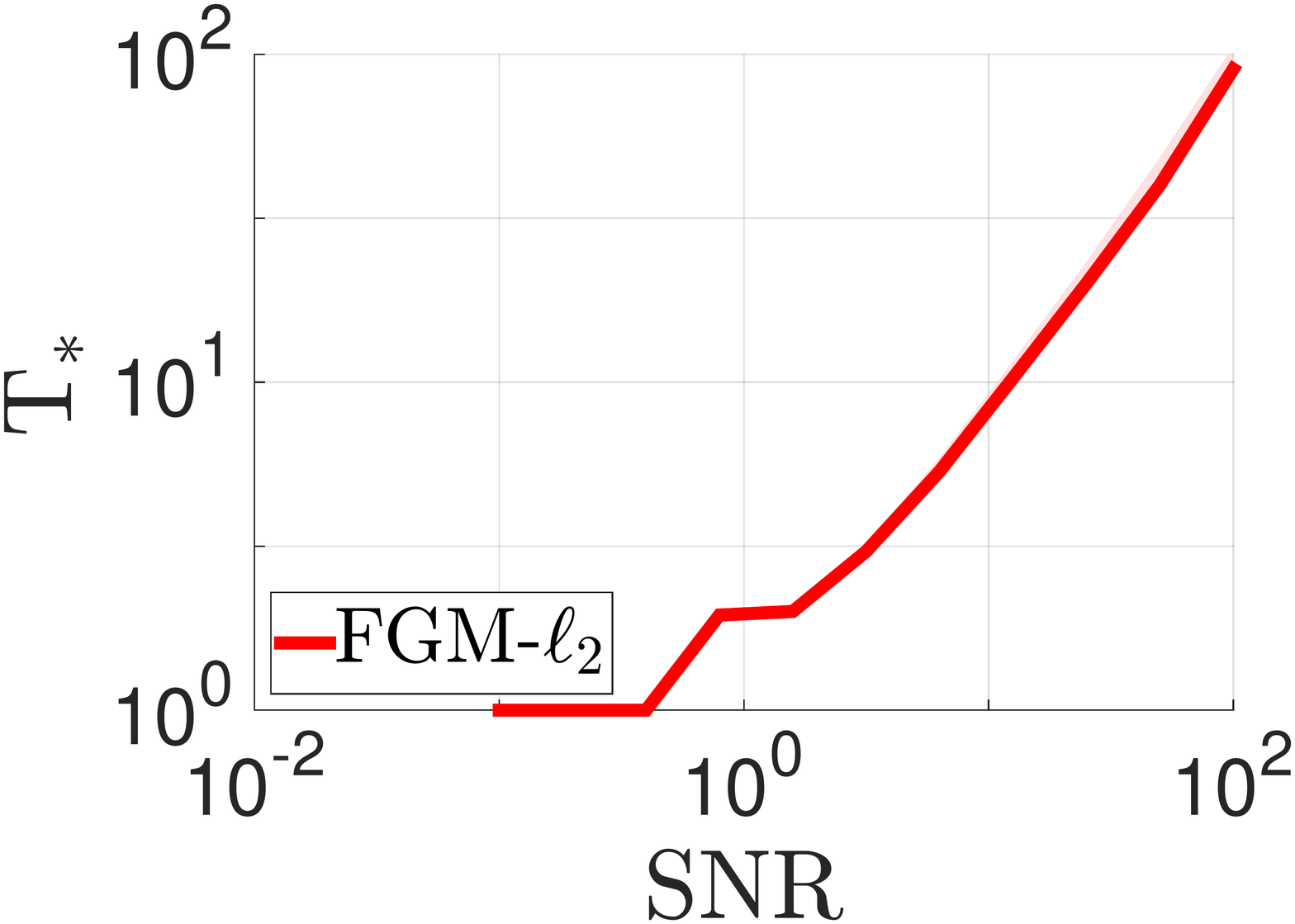}
\end{minipage}
\caption{Iteration at which the accuracy $\veps_*$ is attained for~\eqref{opt:l8con}, left, and~\eqref{opt:l2con}, right, in~\RanSin-4.}
\label{fig:algos-exp-complexity}
\end{figure}

In this experiment (see Fig.~\ref{fig:algos-exp-complexity}), we illustrate the affine dependency of the statistical complexity~$T_*$ from $\SNR$ predicted by our theory, see~\eqref{eq:T-psnr} and~\eqref{eq:psnr-sparse}; note that although the signal in \RanSin~is not sparse on the DFT grid, its DFT is likely to have only a few large spikes which would suffice for~\eqref{eq:psnr-sparse}.
For various $\SNR$ values, we generate a signal in scenario~\RanSin-4, and define the first iteration at which $\veps(T)$ crosses level~$\sigma r$ for~\eqref{opt:l8con} solved with Algorithm~\ref{alg:cmp}, and~$\sigma^2 r^2$ for~\eqref{opt:l2con} with Algorithm~\ref{alg:fgm}. 
We see that the log-log curves plateau for low $\SNR$ and have unit tangent for high $\SNR$, confirming our predictions.

\begin{figure*}
\center
\begin{minipage}{\quartwth}
\centering
\includegraphics[width=1\textwidth, height=\quarthgt, clip=true, angle=0]{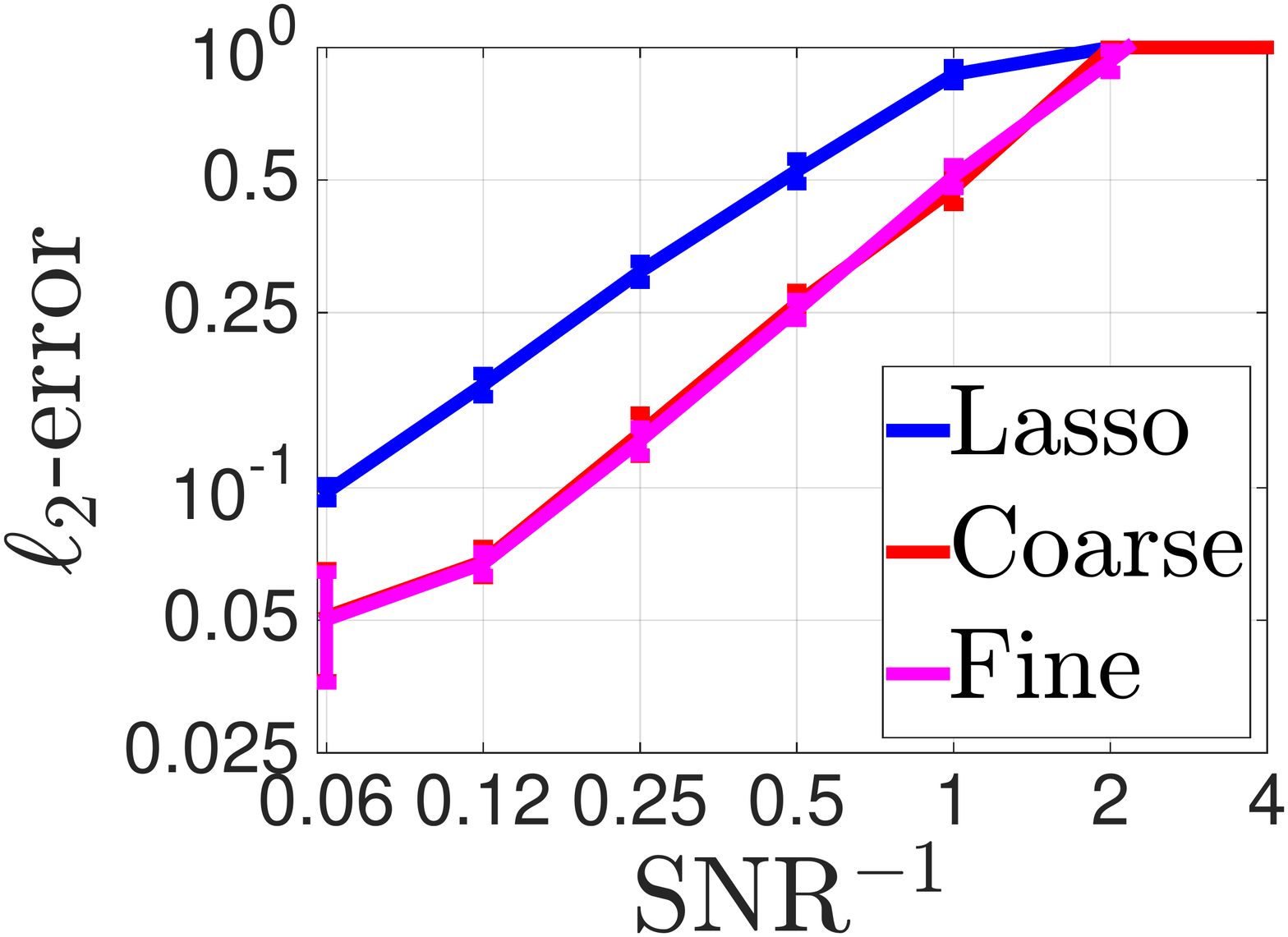}
\includegraphics[width=1\textwidth, height=\quarthgt, clip=true, angle=0]{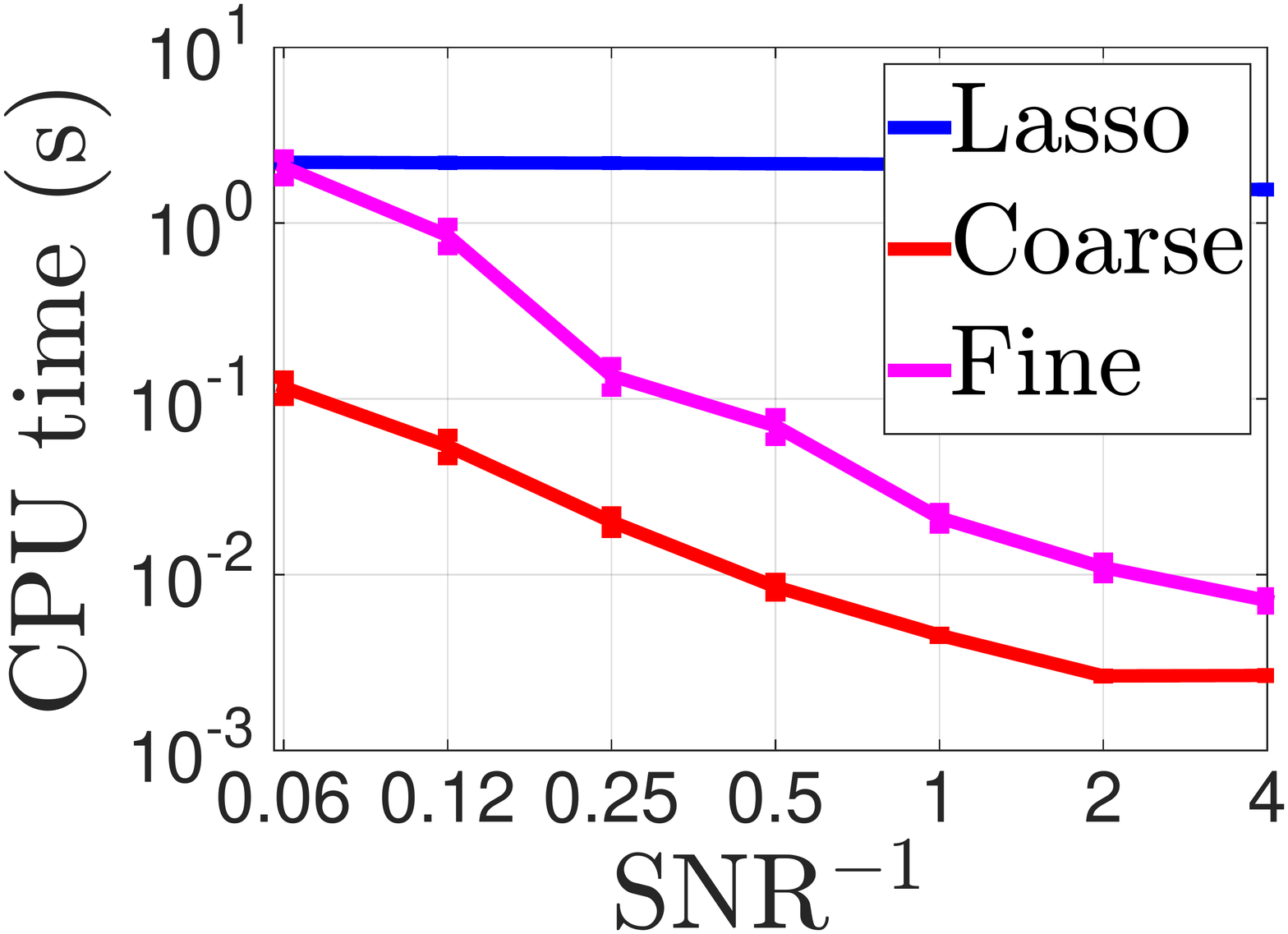}
\end{minipage}
\begin{minipage}{\quartwth}
\centering
\includegraphics[width=1\textwidth, height=\quarthgt, clip=true, angle=0]{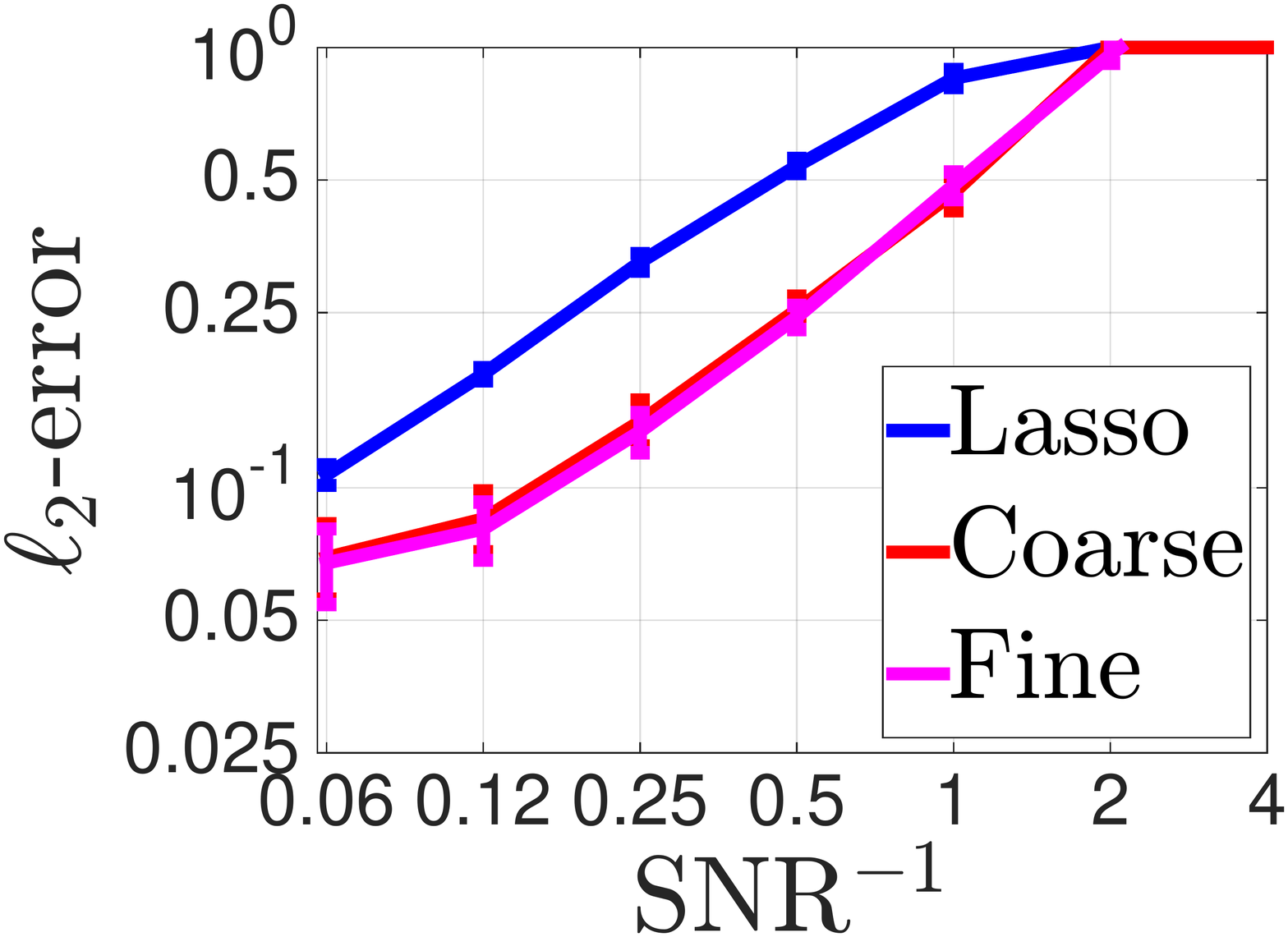}
\includegraphics[width=1\textwidth, height=\quarthgt, clip=true, angle=0]{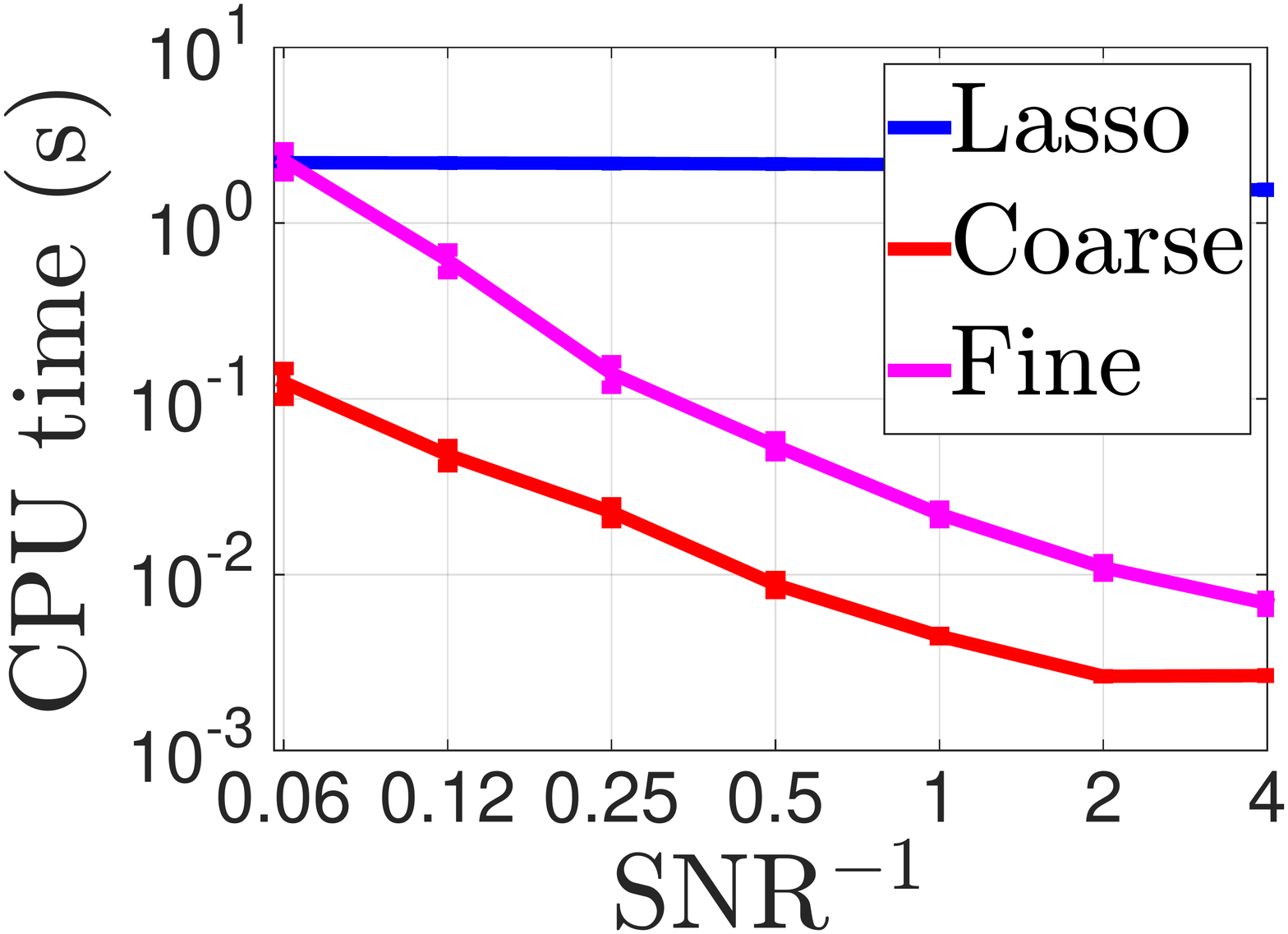}
\end{minipage}
\begin{minipage}{\quartwth}
\centering
\includegraphics[width=1\textwidth, height=\quarthgt, clip=true, angle=0]{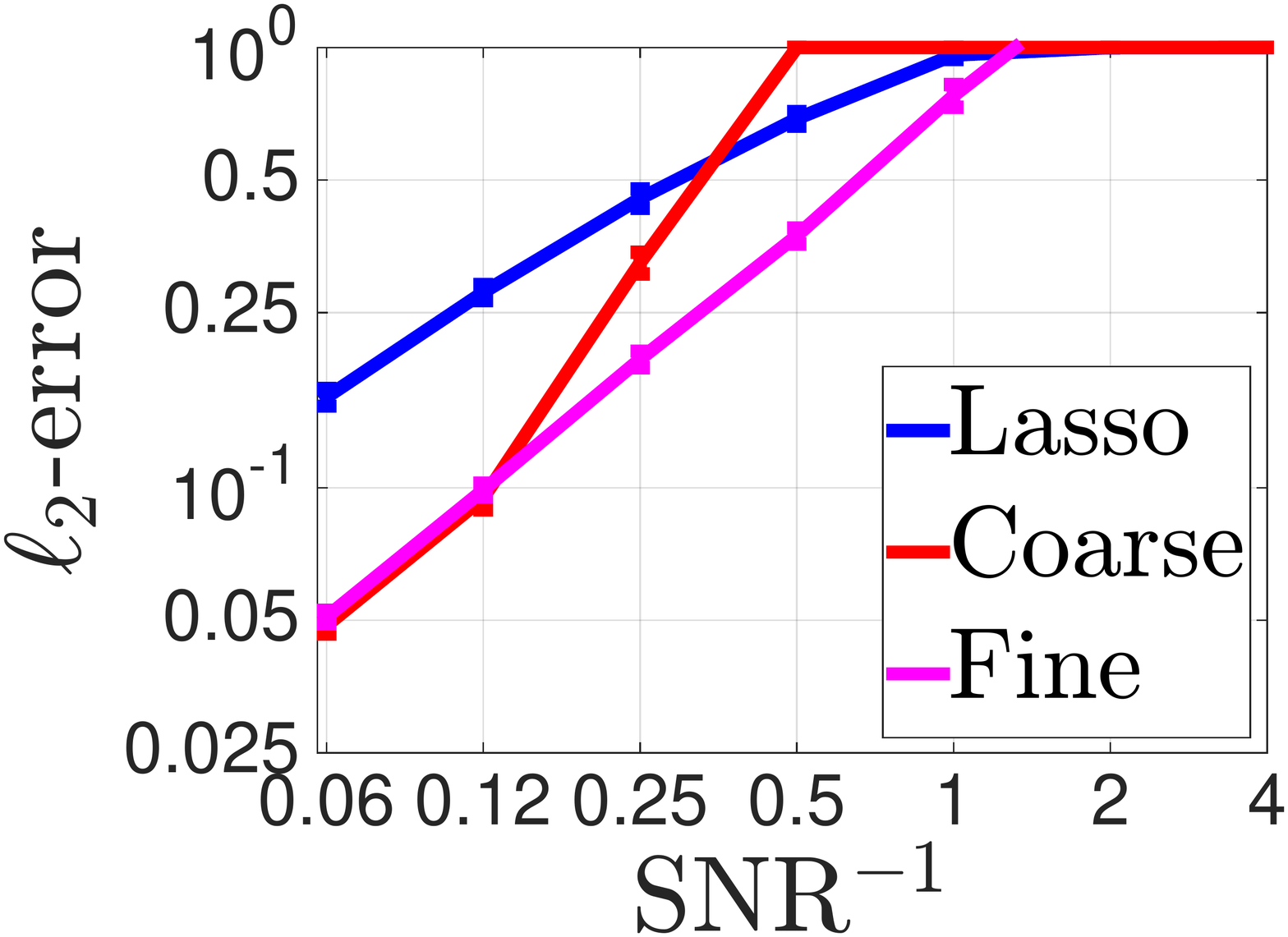}
\includegraphics[width=1\textwidth, height=\quarthgt, clip=true, angle=0]{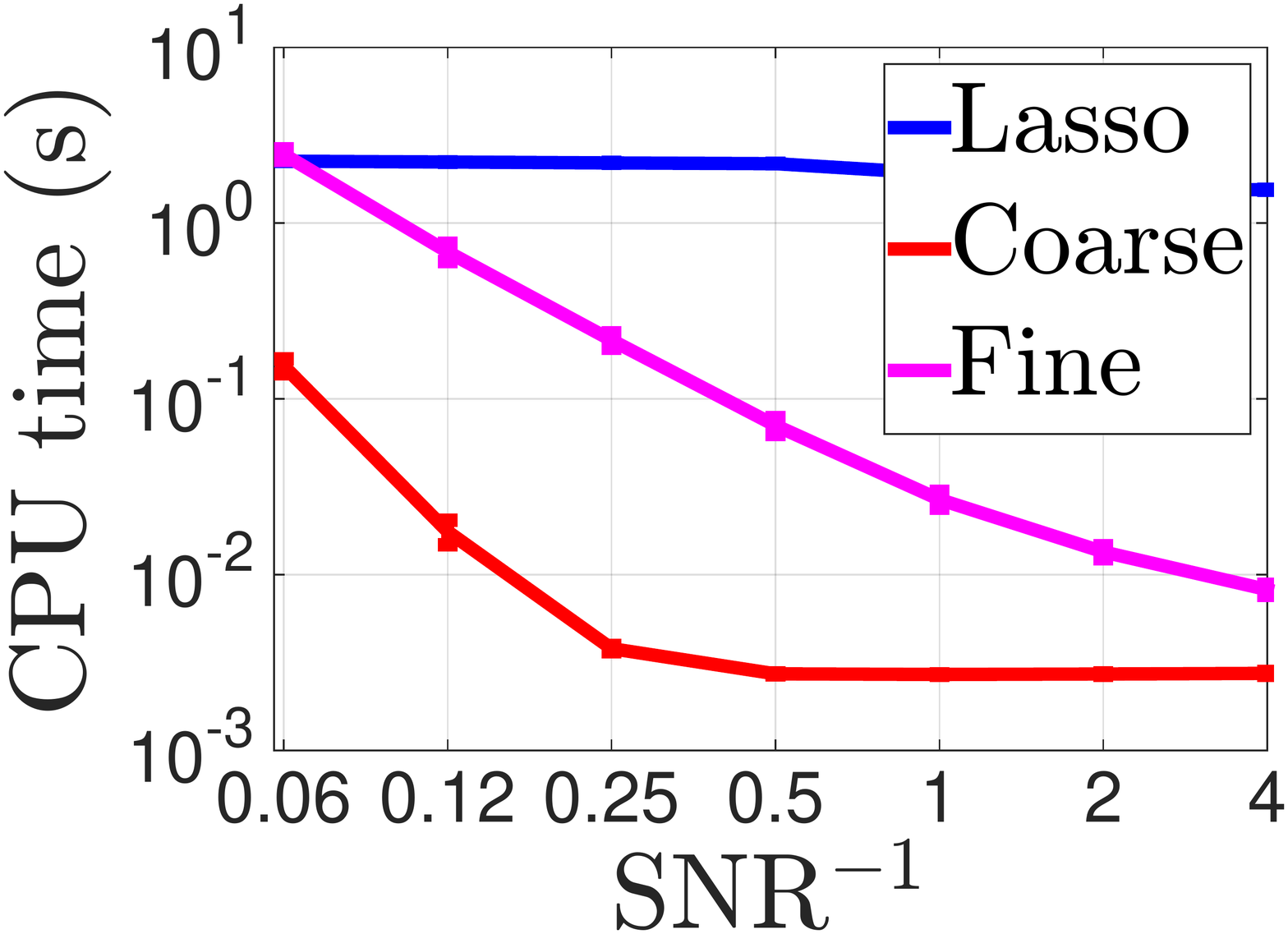}
\end{minipage}
\begin{minipage}{\quartwth}
\centering
\includegraphics[width=1\textwidth, height=\quarthgt, clip=true, angle=0]{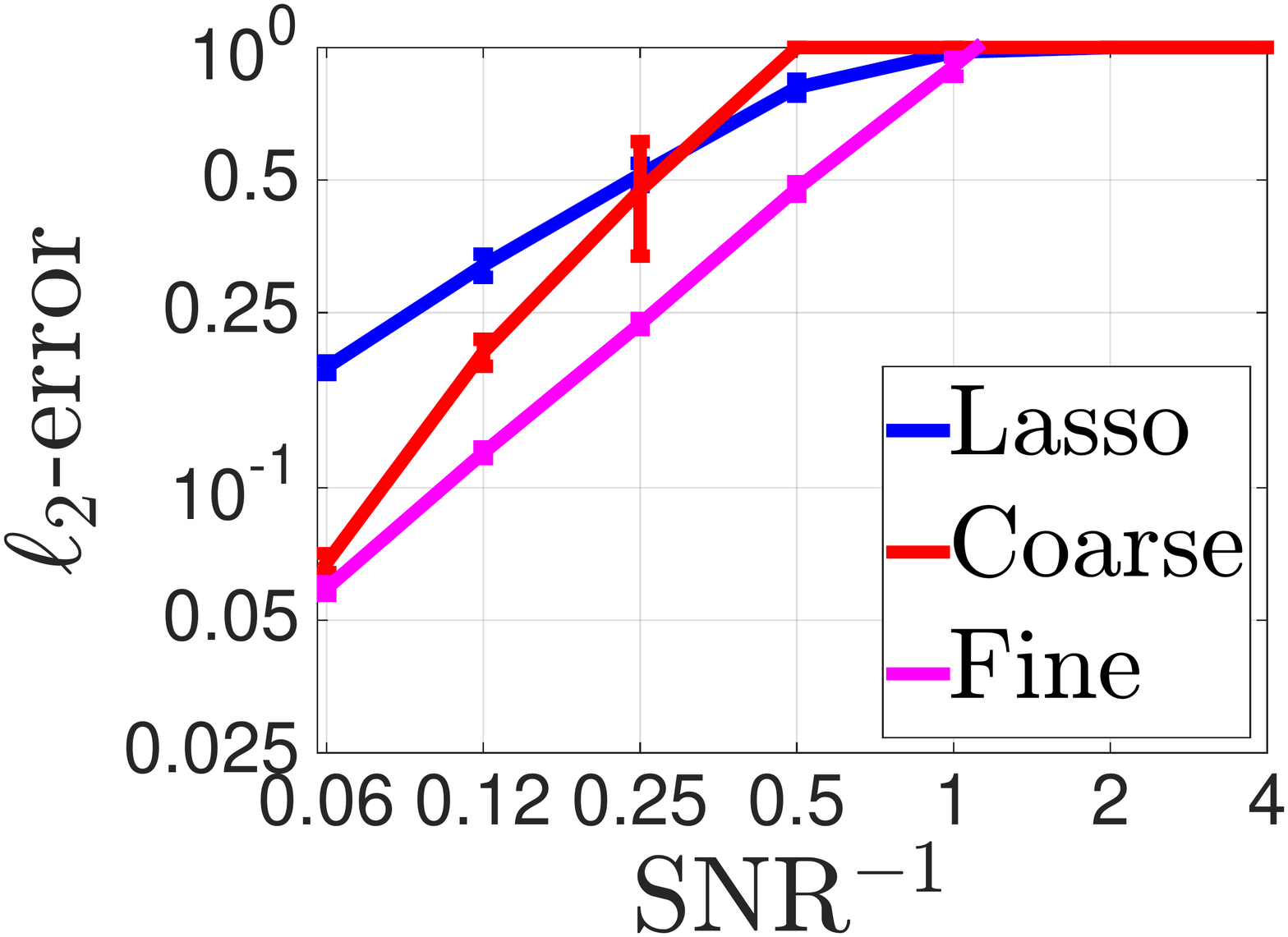}
\includegraphics[width=1\textwidth, height=\quarthgt, clip=true, angle=0]{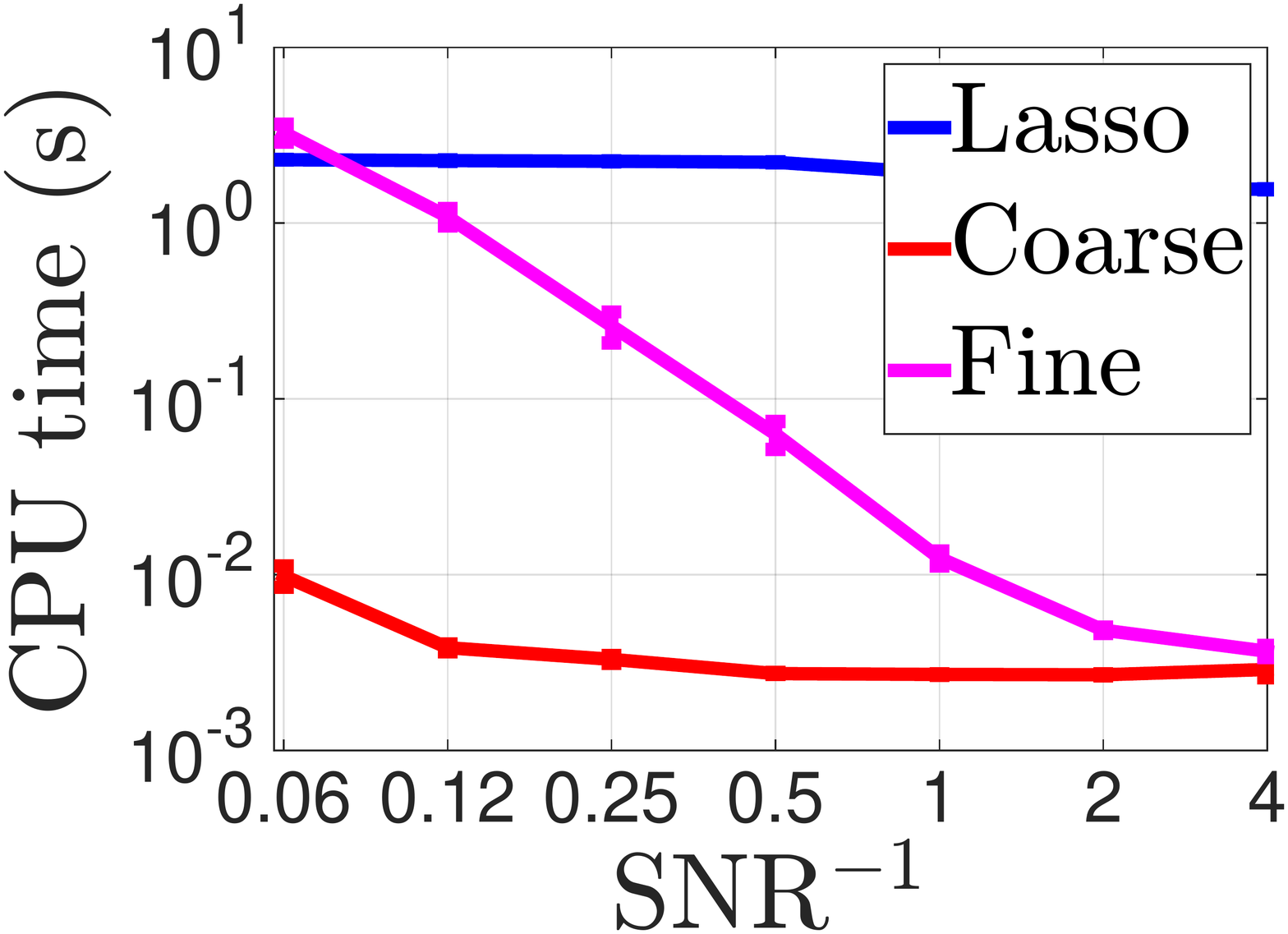}
\end{minipage}

\caption{$\ell_2$-loss and CPU time spent to compute estimators $\vphi^{\coarse}$, $\vphi^{\fine}$, and Lasso.}
\label{fig:algos-exp-sigm}
\end{figure*}

\paragraph{Statistical Performance with Early Stopping.}
In this experiment, we present additional scenario \ModSin-$s$-$m$, in which the signal is a sum of sinusoids with polynomial modulation:
$
x_t = \sum_{k=1}^s p_k(t) e^{i\omega_k t},
$
where $p_k(\cdot)$ are i.i.d. polynomials of degree $r$ with i.i.d. coefficients sampled from~$\C\N(0,1)$; note that in this case~$\dim(\S) = 2s(m+1)$.
Our goal is to study how the early stopping of an algorithm upon reaching accuracy $\veps_*$ (using an accuracy certificate) affects the statistical performance of the resulting estimator. For that, we generate signals in scenarios~\RanSin-4, \CohSin-2, \ModSin-4-2 (quadratic modulation), and \ModSin-4-4 (quartic modulation), with different \SNR, and compare three estimators: approximate solution $\vphi^{\coarse}$ to~\eqref{opt:l2con} with guaranteed accuracy $\veps_*=\sigma^2r^2$, near-optimal solution $\vphi^{\fine}$ with guaranteed accuracy $0.01\veps_*$, and the Lasso estimator, with the standard choice of parameters as described in~\cite{recht1}, which we compute by running 3000 iterations of the FISTA algorithm~\cite{beck2009fast}; 
note that the optimization problem in the latter case is unconstrained, and we do not have an accuracy certificate.
We plot the scaled $\ell_2$-loss of an estimator and the CPU time spent to compute it (we used MacBook Pro 2013 with 2.4 GHz Intel Core i5 CPU and 8GB of RAM).
The results are shown in Fig.~\ref{fig:algos-exp-sigm}.
We observe that $\vphi^{\coarse}$ has almost the same performance as $\vphi^{\fine}$ while being computed 1-2 orders of magnitude faster on average; both significantly outperform Lasso in all scenarios. 

\section*{Acknowledgements}
The authors would like to thank Anatoli Juditsky for fruitful discussions. This work was supported by the LabEx PERSYVAL-Lab (ANR-11-LABX-0025), the project Titan (CNRS-Mastodons), the project MACARON (ANR-14-CE23-0003-01), the NSF TRIPODS Award (CCF-1740551), the program ``Learning in Machines and Brains'' of CIFAR, and a Criteo Faculty Research Award.

\pagebreak
\bibliography{references}
\bibliographystyle{alpha}

\appendix
\section{Background}

\subsection{Adaptive signal denoising}
\label{sec:background-adapt-est}
Assume that the goal is to estimate the signal only on $[0,n]$, from observations~\eqref{eq:intro-observations}, and consider convolution-type estimators
\begin{equation}
\label{eq:background-convolution-est}
\wh x^\vphi_t = [\vphi * y]_t := \sum_{\tau \in \Z} \vphi_\tau y_{t-\tau} \quad 0 \le t \le n.
\end{equation}
Here, $\vphi$ is itself an element of $\C(\Z)$ called a \textit{filter}; 
note that if $\vphi \in \C_n(\Z)$, \eqref{eq:intro-convolution-est} defines an estimator of the projection of $x \in \C(\Z)$ to $\C_n(\Z)$ from observations~\eqref{eq:intro-observations} on $\C_n^{\pm}(\Z)$.
If the filter $\vphi$ is fixed and does not depend on the observations, estimator~\eqref{eq:intro-convolution-est} is linear in observations; otherwise it is not.
Now, assume, following~\cite{ostrovsky2016structure}, that $x \in \C(\Z)$ belongs to a \textit{shift-invariant} linear subspace $\S$ of $\C(\Z)$ 
-- an invariant subspace of the unit shift operator
\[
\Delta : \C(\Z) \to \C(\Z), \quad [\Delta x]_t = x_{t-1}.
\]
As shown in~\cite{ostrovsky2016structure}, one can explicitly construct a filter $\phi^o$, depending on $\S$, such that the worst-case $\ell_2$-risk of the estimator~\eqref{eq:intro-convolution-est} with $\vphi = \phi^o$ satisfies
\begin{equation}
\label{eq:minmaxrate}
\E^{\frac{1}{2}} \left\{\|x - \phi^o * y\|_{n,2}^2 \right\} \le \frac{\sigma\rho}{\sqrt{n+1}} \quad \forall x \in \S,
\end{equation}
where the factor $\rho = \tilde\cO(s^{\kappa})$ for some $\kappa > 0$, that is, is polynomial on the subspace dimension $s = \dim(\S)$ and logarithmic in the sample size (the logarithmic factor can be dropped in some situations).
In fact, one even has a pointwise bound: for any $0\le \tau \le n$, with prob.~$\ge 1-\delta$,
\begin{equation}
\label{eq:minmaxrate-l8}
|x_\tau - [\phi^o * y]_\tau| \le \frac{C\sigma\rho\sqrt{1+\log\left(\frac{n+1}{\delta}\right)}}{\sqrt{n+1}} \quad \forall x \in \S.
\end{equation}
Note that for any fixed subspace $\S$, not even a shift-invariant one, the worst-case $\ell_2$-risk and pointwise risk of \textit{any} estimator can both bounded from below with $c\sqrt{{s}/{(n+1)}}$ for some absolute constant~$c$~\cite{johnstone-book}.
Hence, $\wh x^{\phi^o} = \phi^o * y$ is nearly minimax on $\S$ as long as $s \ll n$: its ``suboptimality factor'' -- the ratio of its worst-case $\ell_2$-risk to that of a minimax estimator -- only depends on the subspace dimension~$s$ but not on the sample size~$n$. 
Unfortunately, $\wh x^{\phi^o}$ depends on subspace $\S$ through the ``oracle'' filter $\phi^o$, and hence it cannot be used in the adaptive estimation setting where the subspace $\S$ with $\dim(\S) = s$ is unknown, but one still would like to attain bounds of the type~\eqref{eq:minmaxrate}.
However, adaptive estimators can be found in the convolution form $\wh x = \wh\vphi * y$ where filter $\wh \vphi = \wh \vphi(y)$ is not fixed anymore, but instead is inferred from the observations. Moreover, $\wh \vphi$ is given as an optimal solution of a certain optimization problem.
Several such problems have been proposed, all resting upon a common principle -- minimization of the Fourier-domain residual 
\begin{equation}
\label{eq:background-residual}
\| F_n[y - \vphi * y]\|_{p}
\end{equation}
with regulzarization via the $\ell_1$-norm $\|F_n[\vphi]\|_{1}$ of the DFT of the filter.
Such regularization is motivated by the following non-trivial fact, see~\cite{harchaoui2015adaptive}: given an oracle filter $\phi^o \in \C_{\lfloor n/2 \rfloor}(\Z)$ which satisfies~\eqref{eq:minmaxrate} with $n$ replaced with $3n$, one can point out a new filter $\vphi^o \in \C_n(\Z)$ which satisfies a ``slightly weaker'' counterpart of~\eqref{eq:minmaxrate-l8},
\begin{equation}
\label{eq:intro-risk}
|x_\tau - [\vphi^o * y]_\tau| \le \frac{3 \sigma r \sqrt{1+\log\left(\frac{n+1}{\delta}\right)}}{\sqrt{n+1}} \quad \forall x \in \S
\end{equation}
where $r = 2\rho^2$, but also admits a bound on DFT in $\ell_1$-norm:
\begin{equation}
\label{eq:intro-constraint}
\|F_n[\vphi^o]\|_{1} \le \frac{r}{\sqrt{n+1}}, \quad r = 2\rho^2.
\end{equation}
see~\cite{ostrovsky2016structure}. 
In fact,~\eqref{eq:intro-constraint} is the key property that allows to control the statistical performance of adaptive convolution-type estimators. 
In some situtaions, polynomial upper bounds on the function $\rho(s)$ are known. 
Then, adaptive convolution-type estimators with provable statistical guarantees can be obtained by minimizing the residual~\eqref{eq:intro-residual} with $p = \infty$~\cite{harchaoui2015adaptive} or $p= 2$~\cite{ostrovsky2016structure} under the constraint~\eqref{eq:intro-constraint}.
A more practical approach is to use penalized estimators, cf. Sec.~\ref{sec:intro}, that attain similar statistical bounds, see~\cite{ostrovsky2016structure} and references therein.

\subsection{Online accuracy certificates}
\label{sec:algos-certificates}
The guarantees on the accuracy of optimization algorithms presented in Section~\ref{sec:algos-general} have a common shortcoming. They are ``offline'' and worst-case, stated once and for all, for the worst possible problem instance. Neither do they get improved in the course of computation, nor become more optimistic when facing an ``easy'' problem instance of the class. 
However, in some situations, online and ``opportunistic'' bounds on the accuracy are available. Following the terminology introduced in~\cite{nemirovski2010accuracy}, such bounds are called \textit{accuracy certificates}. 
They can be used for early stopping of the algorithm if the goal is to reach some fixed accuracy $\varepsilon$). 
One situation in which accuracy certificates are available is saddle-point minimization (via a first-order algorithm) in the case where the domains are bounded and admit an efficiently computable~\textit{linear maximization oracle}. The latter means that the optimization problems
$
\max_{u \in U} \langle a,u\rangle, \quad \max_{v \in V} \langle b,v \rangle
$
can be efficiently solved for any $a, b$. An example of such domains is the unit ball of a norm $\|\cdot\|$ for which the dual norm $\|\cdot\|_*$ is efficiently computable. Let us now demonstrate how an accuracy certificate can be computed in this situation (see~\cite{nemirovski2010accuracy,harchaoui2015conditional} for a more detailed exposition).

A \textit{certificate} is simply a sequence $\lambda^t = (\lambda^t_\tau)_{\tau=1}^t$ of positive weights such that $\sum_{\tau = 1}^t \lambda^t_\tau = 1$.
Consider the $\lambda^t$-average of the iterates $z_\tau$ obtained by the algorithm,
\[
z^t = [u^t, v^t] = \sum_{\tau = 1}^t \lambda^t_\tau z_\tau.
\]
A trivial example of certificate corresponds to the constant stepsize, and amounts to simple averaging. However, one might consider other choices of certificate,
for which theoretical complexity bounds are preserved -- for example, it might be practically reasonable to average only the last portion of the iterates, a strategy called ``suffix averaging''~\cite{rakhlin2012making}.
The point is that any certificate implies a non-trivial (and easily computable) upper bound on the accuracy of the corresponding candidate solution $z^t$. Indeed, the duality gap of a composite saddle-point problem can be bounded as follows:
\begin{equation*}
\begin{aligned}
\overline\phi(u^t) - \underline\phi(v^t) 
&= \overline\phi(u^t) - \phi(u^t,v^t) + \phi(u^t,v^t) - \underline\phi(v^t) \\
&= \max_{v \in V} [\phi(u^t, v) - \phi(u^t,v^t)] - \min_{u \in U} [\phi(u,v^t) - \phi(u^t,v^t)]\\
&\le \max_{v \in V} [\phi(u^t, v) - \phi(u^t,v^t)] +\max_{u \in U} [\phi(u^t,v^t) - \phi(u,v^t)].
\end{aligned}
\end{equation*}
Now, using concavity of $f$ in $v$, we have
\begin{align*}
\phi(u^t, v) - \phi(u^t,v^t) = f(u^t, v) - f(u^t,v^t) \le \sum_{\tau=1}^t \lambda_\tau^t \langle  F_v(z_\tau), v^t-v \rangle.
\end{align*}
On the other hand, by convexity of $f$ and $\Psi$ in $u$,
\begin{align*}
&\phi(u^t,v^t) - \phi(u,v^t) = f(u^t,v^t) - f(u,v^t) + \Psi(u^t) - \Psi(u)  \le \sum_{\tau=1}^t \lambda_\tau^t \langle  F_u(z_\tau) + h(u_\tau), u^t - u \rangle
\end{align*}
where $h(u_\tau)$ is a subgradient of $\Psi(\cdot)$ at $u_\tau$.
Combining the above facts, we get that
\begin{equation}
\label{eq:certificate-accuracy}
\begin{aligned}
\overline\phi(u^t) - \underline\phi(v^t)  \le \max_{u \in U} [-F^t_u - h^t]  + \max_{v \in V} [-F^t_v] + \sum_{\tau = 1}^t \lambda^t_\tau \left[ \langle F_u(z_\tau) + h(u_\tau), u^t\rangle  + \langle F_v(z_\tau), v^t \rangle \right],
\end{aligned}
\end{equation}
where 
\[
F^t_u = \sum_{\tau = 1}^t \lambda^t_\tau F_u(z_\tau), \quad F^t_v = \sum_{\tau = 1}^t \lambda^t_\tau F_v(z_\tau), \quad  \text{and} \;\; h^t = \sum_{\tau = 1}^t \lambda^t_\tau h(u_\tau).
\]
Note that the corresponding averages can often be recomputed in linear time in the dimension of the problem, and then upper bound~\eqref{eq:certificate-accuracy} can be efficiently maintained. 
For example, this is the case when $\lambda^t$ corresponds to a fixed sequence  $\gamma_1, \gamma_2, ...,$ 
\[
\lambda^t_\tau = \frac{\gamma_\tau}{\sum_{\tau' \le t} \gamma_{\tau'}}, \quad \tau \le t.
\]

Note also that any bound on the duality gap implies bounds on the \textit{relative} accuracy for the primal and the dual problem provided that $\underline\phi(v^t)$ (and hence the optimal value  $\phi(u^*,v^*)$) is strictly positive (we used this fact in our experiments, see~Sec.~\ref{sec:experiments}). Indeed, let $\veps(t)$ be an upper bound on the duality gap (\eg~such as~\eqref{eq:certificate-accuracy}), and hence also on the primal accuracy:
\[
\overline\phi(u^t) - \phi(u^*,v^*) \le \overline\phi(u^t) - \underline\phi(v^t)  \le \veps(t).
\]
Then, since~$\phi(u^*,v^*) \ge \underline\phi(v^t) > 0$, we arrive at
\[
\frac{\overline\phi(u^t) - \phi(u^*,v^*)}{\phi(u^*,v^*)} \le \frac{\varepsilon(t)}{\underline{\phi}(v^t)}.
\]
A similar bound can be obtained for the relative accuracy of the dual problem.
\section{Technical proofs}
\label{sec:algos-proofs}

\paragraph{Proof of Lemma~\ref{th:op-norm-upper}.}
Note that $\A$ can be expressed as follows, cf.~\eqref{eq:operator-long-exp}:
\begin{equation}
\label{eq:op-complex-rep}
\A  = \sqrt{2n+1} \cdot F^{\vphantom \H}_n P^{\vphantom \H}_n F^\H_{2n} D^{\vphantom \H}_y F^{\vphantom \H}_{2n} P_n^\H F_n^\H.
\end{equation}
By Young's inequality, for any $\psi \in \C^{n+1}$ we get 
\begin{equation*}
\begin{aligned}
\frac{1}{2n+1}\|\A \psi \|_2^2 
&\le \left\| D^{\vphantom \H}_y F^{\vphantom \H}_{2n} P_n^\H F_n^\H \psi \right\|_2^2 \\
&\le \big\|F_{2n}^{\vphantom\H} [y]_{-n}^n \big\|_\infty^2 \left\| F^{\vphantom \H}_{2n} P_n^\H F_n^\H \psi \right\|_2^2 \\
&\le \big\|F_{2n}^{\vphantom\H} [y]_{-n}^n \big\|_\infty^2 \left\| \psi \right\|_2^2,
\end{aligned}
\end{equation*}
where we used that $P_n$ is non-expansive. 
\qed

\paragraph{Proof of Proposition~\ref{th:op-norm-lower}.}
Consider the uniform grid on the unit circle
\[
U_n = \left\{\exp\left(\frac{2\pi i j}{n+1}\right)\right\}_{j = 0}^{n},
\] 
and the twice finer grid
\[
U_N = \left\{\exp\left(\frac{2\pi i j}{N+1}\right)\right\}_{j = 0}^{N}, \quad N = 2n+1.
\]
Note that $U_N$ is the union of $U_n$ and the shifted grid
\[
\tilde U_n = \left\{u\, e^{i\theta}, \; u \in U_n \right\}, \quad \theta = \frac{2\pi}{N+1};
\]
note that $\tilde U_n$ and $U_n$ do not overlap.
One can check that for any $n \in \Z_+$ and $x \in \C_n(\Z)$, the components of $F_n [x]_0^n$ form the set
\[
\left\{\frac{x(\nu)}{\sqrt{n+1}}\right\}_{\nu \in U_n},
\]
where $x(\cdot)$ is the Taylor series corresponding to $x$:
\[
x(\nu) := \sum_{\tau \in \Z} x_\tau \nu^\tau.
\]
Now, let $x$ be as in the premise of the theorem, and let $x^{(n)} \in \C_n(\Z)$ be such that $x^{(n)}_\tau = x_\tau$ if $0 \le \tau \le n$ and $x^{(n)}_\tau = 0$ otherwise. Similarly, let us introduce $x^{(N)}$ as $x$ restricted on $\C_N(\Z)$. Then one can check that for any $\nu \in U_N$,
\begin{equation}
\label{eq:dft-repeat}
x^{(N)}(\nu) = \left\{
\begin{array}{ll}
2x^{(n)}(\nu), \quad &\nu \in U_n,\\
0, \quad &\nu \in \tilde U_n.
\end{array}
\right.
\end{equation}
In particular, this implies that 
\begin{equation}
\label{eq:concatenation-comparison}
\|F_N[x]_0^N\|_{\infty} = \sqrt{2}\|F_n[x]_0^n\|_{\infty}.
\end{equation}
Now, for any $\vphi \in \C_n(\Z)$, let $\phi \in \C_{N}(\Z)$ be its $n+1$-periodic extension, defined by
\[
[\phi]_0^N = [[\vphi]_0^n; [\vphi]_0^n].
\]
One can directly check that for $x$ as in the premise of the theorem, the circular convolution of $[\phi]_0^N$ and $[x]_0^N$ is simply a one-fold repetition of $2[\vphi * x]_0^n$. Hence, using the Fourier diagonalization property together with~\eqref{eq:concatenation-comparison} applied for  $[\vphi * x]_0^n$ instead of $x_0^n$, we obtain
\begin{equation}
\label{eq:conv2circconv}
\sqrt{N+1} \, \| F_N[x] \odot F_N [\phi]  \|_{\infty} = 2\sqrt{2} \, \|F_n [x * \vphi]\|_{\infty}
\end{equation}
where $a \odot b$ is the elementwise product of $a, b \in \C^{n+1}$.

Finally, note that since $\sigma = 0$, and, as such, $x = y$ a.s., for any $\psi \in \C^{n+1}$ one has:
\[
\A \psi = F_n [x * \vphi], \quad \text{where} \;\; \vphi = F_n^\H[\psi] \in \C_n(\Z).
\]
Hence, using~\eqref{eq:conv2circconv} with such $\vphi$, we arrive at
\begin{align}
\| \A\psi \|_\infty 
&= \|F_n[x * \vphi]\|_{\infty} \nn
&= \frac{\sqrt{n+1}}{2} \| F_N[x] \odot F_N [\phi]  \|_{\infty} \tag*{[by~\eqref{eq:conv2circconv}]} \nn
&= \sqrt{n+1} \| F_n[x] \odot \psi  \|_{\infty} \tag*{[by~\eqref{eq:dft-repeat}]}.
\end{align}
The claim now follows by maximizing the right-hand side in $\psi \in \C^{n+1}: \|\psi\|_{1} \le 1$. \qed

\onecolumn
\subsection{Proof of Theorem~\ref{th:stat-acc-l8}}
\label{app:stat-acc-proof-l8}
The proof is reduced to the following observation: in order to satisfy~\eqref{eq:stat-acc-l8}, it suffices for $\tilde\vphi \in \C_n(\Z)$ to satisfy
\[
\|F_n \tilde\vphi\|_{1} = \cO\left(\frac{r}{\sqrt{n+1}}\right), \quad \|F_n [y - y*\tilde\vphi] \|_\infty = \tilde \cO\left(\sigma r\right), \; \text{where} \; r = 2\rho^2.
\]
This is a rather straightforward remark to the proof of Proposition~4 in~\cite{harchaoui2015adaptive}. 
We give here the proof for convenience of the reader, and also consider the case of the penalized estimator.

\paragraph{Preliminaries.}
Let $\Delta$ be the unit lag operator such that $[\Delta x]_t = x_{t-1}$ for $x \in \C(\Z)$. 
Note that for any filter $\vphi \in \C_n(\Z)$, one can write $\vphi * y = \varphi(\Delta) y$ where $\vphi(\Delta)$ is the Taylor polynomial 
corresponding to $\vphi$:
\[
\vphi(\Delta) := \sum_{\tau \in \Z} \vphi_\tau \Delta^\tau =  \sum_{0 \le \tau \le n} \vphi_\tau \Delta^\tau.
\]
Besides, let us introduce the random variable
\[
\Theta_n(\zeta) := \max_{0 \le \tau \le n}\| \Delta^\tau F_n[\zeta]\|_{\infty}.
\]
Note that  $F_n[\zeta]$ is distributed same as $[\zeta]_0^n$ by the unitary invariance of the law $\C\N(0,I_n)$. 
Using this fact, it is straightforward to obtain that with probability at least $1-\delta$,
\begin{equation}
\label{eq:max-noise}
\Theta_n(\zeta) \le \overline{\Theta}_{n} := 4\sqrt{\log \left(\frac{n+1}{\delta}\right)},
\end{equation}
see~\cite{harchaoui2015adaptive}. 

\paragraph{Constrained uniform-fit estimator.}

Let $\wh\vphi$ be an optimal solution to~\eqref{opt:l8con} with $\overline{r} = r$.
We begin with the following decomposition (recall that $\vphi * y = \varphi(\Delta) y$):
\be
|[x - \widehat \varphi(\Delta) y]_n|
&\le& \sigma |[\widehat \varphi(\Delta) \zeta]_n| + |[x- \widehat\varphi(\Delta) x]_n| \nn
&\le& \sigma \|F_n[\widehat \varphi]\|_{1} \|F_n[\zeta]\|_{\infty} + |[x- \widehat\varphi(\Delta) x]_n| \nn
&\le& \frac{\sigma r  \Theta_n(\zeta)}{\sqrt{n+1}} + |[x- \widehat\varphi(\Delta) x]_n|.
\ee{eq:point-first-decomp}
\noindent 
Here, to obtain the second line we used Young's inequality,
and for the last line we used feasibility of $\wh\vphi$ in~\eqref{opt:l8con}.
Now let us bound $|[x- \widehat\varphi(\Delta) x]_n|$:
\bse
|[x- \widehat\varphi(\Delta) x]_n| 
&\le& |[(1- \widehat\varphi(\Delta))(1-\varphi^o(\Delta))x]_{n}| \; + \; |[\varphi^o(\Delta)(1-\widehat\varphi(\Delta))x]_{n}| \\
&\le& (1 + \|\widehat\varphi\|_{1})\, \|[(1 - \varphi^o(\Delta)) x]_0^n\|_{\infty} \; + \; \|F_n[\varphi^o]\|_{1} \| F_n[(1 - \widehat \varphi(\Delta)) x] \|_{\infty}.
\ese
\noindent Discrepancy of the oracle $\vphi^o$ in the time domain can be bounded using~\eqref{eq:intro-risk}: 
\begin{equation}
\label{eq:l8-est-l8-discrep}
\|[(1 - \varphi^o(\Delta)) x]_0^n\|_{\infty} \le \frac{4r\sigma}{\sqrt{n+1}}.
\end{equation}
Indeed, for any $\tau \in \Z$, $[(1 - \varphi^o(\Delta)) x]_\tau = [x - \vphi^o(\Delta)y]_\tau + \sigma[\vphi^o(\Delta) \zeta]_\tau$. 
On the other hand, using that $\vphi^o$ is non-random,
\[
\E |[\vphi^o(\Delta) \zeta]_\tau|^2 = \|\vphi^o\|_{2}^2 = \|F_n[\vphi^o]\|_{2}^2 \le \|F_n[\vphi^o]\|_{1}^2 =  \frac{r^2}{n+1}.
\]
Now, using that due to~\eqref{eq:intro-constraint} oracle $\varphi^o$ is feasible in~\eqref{opt:l8con}, we can bound the Fourier-domain discrepancy of $\wh\vphi$:
\begin{align}
\|F_n[(1 - \widehat \varphi(\Delta)) x]\|_{\infty}
&\le \left\|F_n[(1 - \widehat \varphi(\Delta)) y]\right\|_{\infty} 
+ \sigma \|F_n[(1 - \widehat \varphi(\Delta)) \zeta]\|_{\infty} \nn
&\le  \left\|F_n[(1 - \widehat \varphi(\Delta)) y]\right\|_{\infty}  
+ \sigma (1 + \|\widehat \varphi\|_1) \Theta_n(\zeta) \nn
&\le \left\|F_n[(1 - \varphi^o(\Delta)) y] \right\|_{\infty}
+ \sigma (1 + \|\widehat \varphi\|_1) \Theta_n(\zeta) \nn
&\le \left\|F_n[(1 - \varphi^o(\Delta)) x]\right\|_{\infty} + \sigma (2+\|\varphi^o\|_1 + \|\widehat \varphi\|_1) \Theta_n(\zeta).
\label{eq:l8-est-discrep}
\end{align}
Meanwhile, using~\eqref{eq:l8-est-l8-discrep}, we can bound the Fourier-domain discrepancy of $\vphi^o$:
\begin{align}\label{eq:l8-orc-discrep}
\left\|F_n[(1 - \varphi^o(\Delta)) x] \right\|_{\infty} 
&\le\left\|F_n[(1 - \varphi^o(\Delta)) x] \right\|_{2}  \nn
&=\left\|[(1 - \varphi^o(\Delta)) x]_0^n\right\|_{2} \le 4\sigma r.
\end{align}
Collecting the above, we obtain
\begin{align*}
|[x- \widehat\varphi(\Delta) x]_n|
&\le (1 + \|\widehat\varphi\|_{1}) \frac{4r \sigma}{\sqrt{n+1}} +\sigma \|F_n[\varphi^o]\|_{1} \left\{4r + (2+\|\varphi^o\|_1 + \|\widehat \varphi\|_1) \Theta_n(\zeta) \right\}.
\end{align*}
Note that $\|F_n[\varphi^o]\|_{1}$ is bounded by~\eqref{eq:intro-constraint}. 
It remains to bound $\|\varphi^o\|_1$ and $\|\widehat \varphi\|_1$:
\begin{equation}\label{eq:l8-filt-norm}
\begin{aligned}
\|\varphi^o\|_1
\le \sqrt{n+1} \|\varphi^o\|_2 
\le \sqrt{n+1} \|F_n[\varphi^o]\|_{1}
\le r,
\end{aligned}
\end{equation}
and similarly $\|\wh \varphi\|_1 \le r$.
Hence, we have
\[
\left|[x - \widehat \varphi(\Delta) x]_n \right| \le \frac{\sigma r}{\sqrt{n+1}} \left[4(1 + 2r) + 2(1 + r) \Theta_n(\zeta)\right],
\]
and, using~\eqref{eq:point-first-decomp} and~\eqref{eq:max-noise}, we arrive that with probability $\ge 1-\delta$,
\begin{equation}
\label{eq:proof-l8con-final}
\left|x_n - [\widehat \varphi(\Delta) y]_n \right| \le \frac{C\sigma r^2\sqrt{1+\log\left(\frac{n+1}{\delta}\right)}}{\sqrt{n+1}}.
\end{equation}
It is now straightforward to see why~$\tilde\vphi$, an $O(\sigma r)$-accurate solution to~\eqref{opt:l8con}, also satisfies~\eqref{eq:proof-l8con-final}: the first change in the above argument when replacing $\wh\vphi$ with $\tilde\vphi$ is the additional term $O(\sigma r)$ in~\eqref{eq:l8-est-discrep}. Since all the remaining terms in the right-hand side of~\eqref{eq:l8-est-discrep} were also bounded from above by~$\cO(\sigma r)$, \eqref{eq:proof-l8con-final} is preserved for $\tilde\vphi$ up to a constant factor. 
\qed

\paragraph{Penalized uniform-fit estimator.}

Let now~$\wh\vphi$ be an optimal solution to~\eqref{opt:l8pen}. 
The proof goes along the same lines as in the previous case; however, we must take into account a different condition for oracle feasibility. 
Proceeding as in~\eqref{eq:point-first-decomp} and using~\eqref{eq:l8-est-l8-discrep}, we get
\begin{equation}
\begin{aligned}
&|[x - \widehat \varphi(\Delta) y]_n| \\
&\;\le \sigma \|F_n[\widehat \varphi]\|_{1} \|F_n[\zeta]\|_{\infty}
+ | [(1 - \widehat \varphi(\Delta)) x ]_n|\\
&\;\le \sigma \|F_n[\widehat \varphi]\|_{1} \|F_n[\zeta]\|_{\infty}
+ \|F_n[\varphi^o]\|_{1}\|F_n[(1 - \widehat \varphi(\Delta)) x]\|_{\infty} 
+  (1 + \|\widehat\varphi\|_{1}) \|[(1 - \varphi^o(\Delta)) x]_0^n\|_{\infty}\\
&\;\le \sigma \|F_n[\widehat \varphi]\|_{1} \Theta_n(\zeta)
+ \frac{r}{\sqrt{n+1}} \|F_n[(1 - \widehat \varphi(\Delta)) x]\|_{\infty}
+  \frac{4r \sigma}{\sqrt{n+1}}(1 + \|\widehat\varphi\|_{1}).
\end{aligned}
\label{eq:l8pen-init}
\end{equation}
Let us condition on the event $\Theta_n(\zeta) \le \overline\Theta_{n}$ the probability of which is $\ge 1-\delta$.
Feasibility of $\wh\varphi$ in~\eqref{opt:l8pen} yields
\begin{align}
\left\| F_n[(1 - \wh\varphi(\Delta)) y]\right\|_{\infty}
+ \lambda \|F_n[\wh\varphi]\|_{1} 
&\le \left\| F_n[(1 - \varphi^o(\Delta)) y] \right\|_{\infty} + \lambda \|F_n[\varphi^o]\|_{1} \nn
&\le 4 \sigma r + (1+r) \sigma \Theta_n(\zeta) + \frac{\lambda r}{\sqrt{n+1}} \nn
&\le \left(4 + 2\Theta_n(\zeta) + \frac{\lambda}{\sigma \sqrt{n+1}}\right) \sigma r \nn
&\le \frac{2\lambda r}{\sqrt{n+1}}.
\label{eq:l8pen-long}
\end{align} 
Here first we used~\eqref{eq:l8-orc-discrep},~\eqref{eq:l8-filt-norm}, and the last line of~\eqref{eq:l8-est-discrep}, then that $r \ge 1$, and, finally, used the choice of $\lambda$ from the premise of the theorem.
Now from~\eqref{eq:l8pen-long} we obtain
\begin{align}
\label{eq:l8pen-norm-phi}
\|F_n[\wh\varphi]\|_{1} \le \frac{2r}{\sqrt{n+1}}
\end{align}
and
\begin{align}
\label{eq:l8pen-norm-aphi}
1+\|\widehat \varphi\|_1 
&\le 1+ \sqrt{n+1} \|F_n[\widehat \varphi]\|_{1}
\le 1 + 2r \le 3r.
\end{align}
Further, using~\eqref{eq:l8pen-long}~and~\eqref{eq:l8pen-norm-aphi}, we get
\begin{align}
\| F_n[(1 - \widehat \varphi(\Delta)) x] \|_{\infty} 
&\le \left\| F_n[(1 - \widehat \varphi(\Delta)) y] \right\|_{\infty}  
+ \sigma (1 + \|\widehat \varphi\|_1) \Theta_n(\zeta)\nn
&\le \left(\frac{2\lambda}{\sigma\sqrt{n+1}} + 3\Theta_n(\zeta)\right) \sigma r
\label{eq:l8pen-last}
\end{align}
Substituting~\eqref{eq:l8pen-norm-phi}--\eqref{eq:l8pen-last} into~\eqref{eq:l8pen-init}, we arrive at
\begin{align*}
|[x - \widehat \varphi(\Delta) y]_0| 
&\le \left(\frac{2\lambda}{\sigma\sqrt{n+1}} + 5\Theta_n(\zeta) + 8\right) \frac{\sigma r^2}{\sqrt{n+1}} \\
&\le \frac{5\lambda r^2}{n+1} \\
&= \frac{80 r^2 \sqrt{1 + \log\left(\frac{n+1}{\delta}\right)}}{\sqrt{n+1}}.
\end{align*}
Similarly to the case of the constrained estimator, it is straightforward to see that the last bound is preserved (up to a constant factor) for an $\veps$-accurate solution~$\tilde\vphi$ to~\eqref{opt:l8pen} with $\veps = O(\sigma r)$. 
\qed

\subsection{Proof of Theorem~\ref{th:stat-acc-l2}}
\label{app:stat-acc-proof-l2}

\paragraph{Constrained least-squares estimator.}
Let us first summarize the original proof of~\eqref{eq:stat-acc-l2} for the case of an exact optimal solution~$\wh\vphi$ of~\eqref{opt:l2con}, see Theorem 2.2 in~\cite{ostrovsky2016structure} and its full version~\cite{arxiv2016structure}.
Introducing the scaled Hermitian dot product for $\vphi, \psi \in \C_n(\Z)$, 
\[
\langle \vphi, \psi \rangle_n = \frac{1}{n+1} \sum_{\tau = 0}^n \overline\vphi_\tau \psi_\tau,
\]
 the squared $\ell_2$-loss can be decomposed as follows:
\begin{align}
\label{eq:stat-acc-chain-0}
\|x-\wh\vphi*y\|_{n,2}^2
&= \|y-\wh\vphi*y\|_{n,2}^2 - \sigma^2 \|\zeta\|_{n,2}^2 - 2\sigma\langle \zeta, x - \wh\vphi * y\rangle_n \nn
&\le \|y-\vphi^o*y\|_{n,2}^2 - \sigma^2 \|\zeta\|_{n,2}^2 - 2\sigma\langle \zeta, x - \wh\vphi * y\rangle_n \nn
&= \|x-\vphi^o*y\|_{n,2}^2 + 2\sigma\langle \zeta, x - \vphi^o * y\rangle_n - 2\sigma\langle \zeta, x - \wh\vphi * y\rangle_n,
\end{align}
where the inequality is due to feasibility of $\vphi^o$ in~\eqref{opt:l2con}. 
Now, it turns out that the dominating term in the right-hand side is the first one (corresponding to the squared oracle loss): we know that due to~\eqref{eq:intro-risk}, with probability $\ge 1-\delta$ one has
\begin{equation}
\label{eq:stat-acc-l2-oracle-risk}
\|x-\vphi^o*y\|_{n,2}^2 \le \frac{9\sigma^2 r^2 \log\left(\frac{n+1}{\delta}\right)}{n+1}.
\end{equation}
On the other hand, one can bound the next term in the right-hand side of~\eqref{eq:stat-acc-chain-0} as 
\begin{align}
\sigma \langle \zeta, x - \vphi^o * y\rangle_n 
&\le \frac{\sigma \sqrt{2 \log\left(\frac{3}{\delta}\right)}}{\sqrt{n+1}} \|x-\vphi^o*y\|_{n,2} + \frac{12\sigma^2 r(1 + \log\left(\frac{6}{\delta}\right))}{n+1} \nn
&\le \frac{6\sigma^2 r \log\left(\frac{3(n+1)}{\delta}\right)}{n+1} + \frac{12\sigma^2 r(1 + \log\left(\frac{6}{\delta}\right))}{n+1} \nn
&\le \frac{30\sigma^2 r \log\left(\frac{6(n+1)}{\delta}\right)}{n+1}.
\label{eq:stat-acc-l2-oracle-cross}
\end{align}
Here, for the first inequality we refer the reader to the original proof in~\cite{arxiv2016structure}, eq.~(44-45), where one should set $\kappa_{m,n} = 1$ and keep in mind the absence of scaling factor $\frac{1}{n+1}$ in the definitions of~ $\langle \phi, \psi \rangle_n$ and $\|\cdot\|_{n,2}$. The next inequalities then follow by simple algebra using~\eqref{eq:stat-acc-l2-oracle-risk}. 

Finally, the last term in the right-hand side of~\eqref{eq:stat-acc-chain-0} can be bounded as follows with probability $\ge 1-\delta$:
\begin{equation}
\label{eq:stat-acc-l2-est-cross}
2\sigma |\langle \zeta, x - \wh\vphi * y\rangle_n|
\le \frac{2\sqrt{2}\sigma\left(\sqrt{s} + \sqrt{\log\left(\frac{2}{\delta}\right)}\right)}{\sqrt{n+1}} \|x - \wh\vphi * y \|_{n,2} 
+ \frac{8\sqrt{2}\sigma^2r \left(2 + \log\left(\frac{8(n+1)}{\delta} \right) \right)}{n+1}.
\end{equation}
see eq.~(33-40) in~\cite{arxiv2016structure} where one must set $\varkappa = 0$ in our setting since $x \in \S$. 
Moreover, in the proof of~\eqref{eq:stat-acc-l2-est-cross} the optimality of $\wh\vphi$ was not used; instead, the argument in~\cite{arxiv2016structure} relied only on the following facts:
\begin{itemize}
\item[(i)] $x \in \S$ where $\S$ is a shift-invariant subspace of $\C(\Z)$ with $\dim(\S) = s$; 
\item[(ii)] one has a bound on the Fourier-domain $\ell_1$-norm of $\widehat\vphi$:
$
\|F_n[\wh\vphi]\|_{1} \le \frac{r}{\sqrt{n+1}}.
$
\end{itemize}
Finally, collecting \eqref{eq:stat-acc-chain-0}-\eqref{eq:stat-acc-l2-est-cross} and solving the resulting quadratic inequality, one bounds the scaled $\ell_2$-loss of $\wh\vphi$:
\begin{equation}
\label{eq:stat-acc-l2-final}
\|x-\wh\vphi*y\|_{n,2} \le \frac{C\sigma}{\sqrt{n+1}} \left( \sqrt{s}+ r \sqrt{\log \left(\frac{n+1}{\delta}\right)} \right).
\end{equation}
(We used that $r \ge 1$.) 
Moreover, it is now evident that an $\veps$-accurate solution~$\wh\vphi$ to~\eqref{opt:l2con} with $\veps = \cO(\sigma^2r^2)$ still satisfies~\eqref{eq:stat-acc-l2-final}. 
Indeed, the error decomposition~\eqref{eq:stat-acc-chain-0} must now be replaced with
\begin{equation}
\label{eq:stat-acc-chain-1}
\|x-\tilde\vphi*y\|_{n,2}^2 
\le \|x-\vphi^o*y\|_{n,2}^2 + 2\sigma\langle \zeta, x - \vphi^o * y\rangle_n - 2\sigma\langle \zeta, x - \tilde\vphi * y\rangle_n + \frac{\veps}{n+1}.
\end{equation}
Then, \eqref{eq:stat-acc-l2-oracle-risk}~and~\eqref{eq:stat-acc-l2-oracle-cross} do not depend on~$\tilde\vphi$, and hence are preserved. 
The term $\frac{\veps}{n+1}$ enters additively, and allows for the same upper bound as~\eqref{eq:stat-acc-l2-oracle-risk}. 
Finally, ~\eqref{eq:stat-acc-l2-est-cross} is preserved when replacing~$\wh\vphi$ with~$\tilde\vphi$ since (i) and (ii) remain true. \qed

\paragraph{Penalized least-squares estimator.}
Let now $\tilde \vphi$ be an $\veps$-accurate solutions to~\eqref{opt:l2pen}, let $\lambda_n = \frac{\lambda}{\sqrt{n+1}}$, and let $\tilde r = \sqrt{n+1}\|F_n[\tilde\vphi]\|_1$.
Similarly to~\eqref{eq:stat-acc-chain-1}, one has
\begin{align}
\label{eq:stat-acc-chain-l2pen}
\|x-\tilde\vphi*y\|_{n,2}^2 \le \|x-\vphi^o*y\|_{n,2}^2 + 2\sigma\langle \zeta, x - \vphi^o * y\rangle_n - 2\sigma\langle \zeta, x - \tilde\vphi * y\rangle_n
  +  \frac{\lambda_n (r - \tilde r)}{n+1} + \frac{\veps}{n+1}.
\end{align}
Note that~\eqref{eq:stat-acc-l2-oracle-risk} and~\eqref{eq:stat-acc-l2-oracle-cross} are still valid. 
Moreover, \eqref{eq:stat-acc-l2-est-cross} is preserved for $\tilde\vphi$ if $r$ is replaced with $\tilde r$, cf.~(i) and~(ii):
\begin{equation}
\label{eq:stat-acc-l2-est-cross-pen}
2\sigma |\langle \zeta, x - \tilde\vphi * y\rangle_n|
\le \frac{2\sqrt{2}\sigma\left(\sqrt{s} + \sqrt{\log\left(\frac{2}{\delta}\right)}\right)}{\sqrt{n+1}} \|x - \tilde\vphi * y \|_{n,2} 
+ \frac{8\sqrt{2}\sigma^2 \tilde r \left(2 + \log\left(\frac{8(n+1)}{\delta} \right) \right)}{n+1}.
\end{equation}
Hence, if $\lambda$ is chosen as in the premise of the theorem, the second term in the right-hand side is dominated by $\frac{\lambda_n \tilde r}{n+1}$. 
Combining \eqref{eq:stat-acc-l2-oracle-risk}, \eqref{eq:stat-acc-l2-oracle-cross}, and~\eqref{eq:stat-acc-l2-est-cross-pen} with the fact that $\veps = \cO(\sigma^2r^2)$, plugging in the value of $\lambda$ from the premise of the theorem, and solving the resulting quadratic inequality, we conclude that~\eqref{eq:stat-acc-l2-final} is preserved for $\tilde\vphi$. \qed

\end{document}